\documentclass{amsart}
\usepackage[headings]{fullpage}
\usepackage{amssymb,latexsym}
\usepackage{amsmath,amscd, rotating}
\usepackage{graphics}
\usepackage{subfigure}
\usepackage[all, knot,line]{xy}

\newtheorem{thm}{Theorem}
\newtheorem{lem}[thm]{Lemma}
\newtheorem{cor}[thm]{Corollary}
\newtheorem{prop}[thm]{Proposition}

\theoremstyle{plain}

\theoremstyle{definition}

\newtheorem{defi}[thm]{Definition}

\newtheorem{problem}{Problem}

\theoremstyle{remark}

\renewcommand{\int}{\operatorname{int}}

\newcommand{\knot}{{k}}
\newcommand{\unknot}{\text{unknot}}
\newcommand{\mfld}{M}
\newcommand{\move}{\mu}
\newcommand{\sphere}{{S^3}}
\newcommand{\m}{\mathcal G}
\newcommand{\mugroup}{G^\move}
\newcommand{\double}{\eta}
\newcommand{\ws}{\mathcal H}
\newcommand{\gv}[1]{{V}_{#1}}

\newcommand{\group}[2]{G^{#1}_{#2}}

\newcommand{\A}{\mathcal A}
\newcommand{\zlog}{\log_{\Z}}

\newcommand{\Z}{\mathbb{Z}}

\newcommand{\Q}{\mathbb{Q}}
\newcommand{\F}{\mathcal{F}}

\newcommand{\subgroup}{\m}

\newcommand{\B}{\mathcal B}

\newcommand{\V}[1]{\F^{V}_{#1}}
\newcommand{\K}{\text{Knots}}
\newcommand{\BP}{\operatorname{BP}}

\newcommand{\DDelta}{\Delta\Delta}
\newcommand{\Null}{\mathcal N}

\newcommand{\s}{{\mathcal S}}
\newcommand{\p}[1]{{\mathcal P}_{#1}}
\newcommand{\one}{\boldsymbol{1}}

\begin{document}
\title[Finite-type invariants based on local moves]{Finite-type knot invariants based on the band-pass and doubled-delta moves}
\author[J. Conant]{James Conant}
\address{Department of Mathematics\\
University of Tennessee at Knoxville\\
Knoxville, TN , 37996}
\email{jconant@math.utk.edu}

\author[J. Mostovoy]{Jacob Mostovoy}
\address{Departamento de Matem\'aticas\\
 CINVESTAV, Apartado Postal 14-740\\
C.P.\ 07000 M\'exico, D.F. Mexico
}
\email{jacob@math.cinvestav.mx}

\author[T. Stanford]{Ted Stanford}
\address{Department of Mathematics\\
New Mexico State University\\
Las Cruces, NM 88003}
\email{stanford@nmsu.edu}


\thanks{The first author was partially supported by NSF grant
DMS 0305012.}

\begin{abstract}
We study generalizations of finite-type knot invariants obtained by
replacing the crossing change in the Vassiliev skein relation by
some other local move.

First, we represent the local moves by normal subgroups of the pure
braid group $\p{\infty}$. Subgroups that are stable under the
``strand-tripling'' endomorphisms are shown to produce finite-type
invariants with familiar properties; in particular, generalized
Goussarov's $n$-equivalence classes of knots form groups under the
connected sum. (Similar results, but with a different approach, have
been obtained before by Taniyama and Yasuhara.)

Treating local moves as surgeries on claspers, we study two
particular cases in detail: the band-pass and the doubled delta
move. While the band-pass move gives only one ``new'' invariant
(namely, the Arf invariant), the invariants corresponding to the
doubled-delta move contain information which is not available to any
finite collection of Vassiliev invariants.

The complete degree 0 doubled-delta invariant is the $S$-equivalence
class of the knot. In this context, we generalize a result of
Murakami and Ohtsuki to show that the only primitive Vassiliev
invariants of $S$-equivalence taking values in an abelian group with
no 2-torsion arise from the Alexander-Conway polynomial. To this
end, we introduce a discrete logarithm which transforms the
coefficients of the Conway polynomial into primitive integer-valued
invariants. As for the higher degree doubled-delta invariants, we
start analyzing them by considering which Vassiliev invariants are
degree 1 as doubled-delta invariants. We find that there is exactly
one Vassiliev invariant in each odd degree which is of doubled-delta
degree one, whereas in each even degree there is at most a
$\Z_2$-valued invariant, which we show exists in degree 4. For
higher doubled delta degrees, we observe that the Euler degree $n+1$
part of Garoufalidis and Kricker's rational lift of the Kontsevich
integral is a doubled-delta degree $2n$ invariant. Finally, the
doubled-delta move is a special case of Garoufalidis and Rozansky's
null-move of pairs $(M; k)$ where $M$ is a homology 3-sphere and $k$
a knot in $M$. A consequence of our work is that n-equivalence
classes of pairs $(S^3; k)$ with respect to the null-move, do, in
fact, form a group.
\end{abstract}

\maketitle

   %%%%%%%%%%%%%%%%%%%%%%%%%%%%%%%%%%%%%%%%%%%%%%%%%%%%%%%
   % uno.                                                %
   %                          %%%                        %
   %                         % %%                        %
   %                        %  %%                        %
   %                       %   %%                        %
   %                           %%                        %
   %                           %%                        %
   %                           %%                        %
   %                           %%                        %
   %                           %%                        %
   %                       %%%%%%%%%                     %
   %                                                     %
   %%%%%%%%%%%%%%%%%%%%%%%%%%%%%%%%%%%%%%%%%%%%%%%%%%%%%%%

\section{Introduction}

It has been observed that the theory of finite-type knot invariants
can be generalized by replacing the crossing change appearing in
the Vassiliev skein relation by some other local move.

Vassiliev's knot-space definition of finite-type invariants provides
little motivation for such generalization. Indeed, a crossing change
on a knot diagram corresponds to a codimension one singularity of a map from
a circle to a 3-sphere, while other local moves correspond to
singularities of higher codimension.

Nevertheless, the theories of finite-type invariants based on local
moves other than crossing changes fit perfectly with the approach
developed in the works of Goussarov \cite{gouss, gouss2, gouss3}
and later, Habiro \cite{habiro} and Stanford \cite{stanford}. In
particular, Goussarov's notion of $n$-equivalence can be extended to
a wide class of local moves; this generalized notion of
$n$-equivalence is consistent with suitably defined finite-type
invariants, and the $n$-equivalence classes of knots form groups
under the connected sum operation. See Stanford \cite{st-delta} for a
particular example of such situation and Taniyama and Yasuhara
\cite{ty2} for a more general treatment of theories based on local
moves.

The local moves considered in  \cite{st-delta} and \cite{ty2}
essentially consist of replacing one fixed subtangle of a link by
another fixed subtangle. A different approach to local moves can be
taken by defining local moves as surgeries on claspers. These local
moves can change the ambient manifold. Finite-type invariants based
on a move of this type (the ``null-move") were considered by
Garoufalidis and Rozansky in \cite{gr}. It was shown in \cite{gr, gk}
that the rational lift of the Kontsevich integral is a
universal rational finite-type (with respect to the null-move)
invariant of knots in integral homology spheres with trivial Alexander polynomial.

Apart form the null-move, the only local moves for which the
corresponding finite-type invariants have been identified to some
extent are the $C_k$-moves. Taniyama and Yasuhara proved in
\cite{ty1} that the space of primitive $C_k$-finite-type invariants
of order $n$ and smaller coincides with the space of usual primitive
finite-type invariants of order $kn$ and smaller.

\medskip

The present paper consists of two parts. In the first part we show
how local moves on knots can be interpreted via pure braid closures.
The notion of a local move is replaced here by the notion of a
modification by an element of a normal subgroup of the pure braid
group on an infinite number of strands $\p{\infty}$ (This is the direct limit
of finite pure braid groups, each included in the next by adding a trivial strand.) 
Each normal
subgroup $\subgroup\subset\p{\infty}$ gives rise to
$\subgroup$-finite-type knot invariants and to the relation of
$\gamma_n\subgroup$-equivalence on the isotopy classes of knots. (Here $\gamma_n$ refers to the nth term of the lower central series.)
Using a Markov-type theorem for the short-circuit closure of
\cite{ms}, we show the following
\begin{thm}\label{thm:main0}
If the subgroup $\subgroup$ is stable under ``strand-tripling", then
\begin{enumerate}
\item two $\subgroup$-trivial knots are $\gamma_n\subgroup$-equivalent
if and only if they cannot be distinguished by
$\subgroup$-finite-type invariants of orders $n$ and smaller;
\item  the set of $\subgroup$-trivial knots modulo
$\gamma_n\subgroup$-equivalence is an abelian group.
\end{enumerate}
\end{thm}
Here a knot is said to be $\subgroup$-trivial if it is
$\subgroup$-equivalent to the unknot. The precise definition of the
strand-tripling endomorphisms of $\p{\infty}$ is given in
Section~\ref{section:moves}.

The above result generalizes the well-known theorems of Goussarov
\cite{gouss}. It is very similar, though not immediately equivalent,
to a theorem by Taniyama and Yasuhara \cite{ty2}. It is probable
that our result can be obtained using methods of \cite{ty2};
nevertheless, we think that our approach is of independent interest.

We also re-state Theorem~\ref{thm:main0} in terms of claspers; this
enables us to apply it to finite-type invariants based on the null
move.

\medskip

In the second part we study the finite-type invariants based on two
particular local moves: the band-pass move and the doubled-delta
move.

\medskip

We prove that the finite-type invariants with respect to the band-pass move essentially
coincide with the Vassiliev invariants. More precisely, primitive
finite-type invariants with respect to the band-pass move of order $n$ coincide with
primitive Vassiliev invariants of order $n$ for $n\geq 1 $. The Arf
invariant is the unique band-pass-finite-type invariant of order
$0$.

\medskip

The case of the doubled-delta move  turns out to be more
interesting.

In what follows, we shall abbreviate ``doubled-delta-finite-type
invariants of order $n$" to ``$\DDelta_n$-invariants". These are
quite interesting already for $n=0$: it has been proved by Naik and
Stanford in \cite{ns} that two knots are $S$-equivalent if and only
if one can be transformed to the other by a sequence of doubled
delta moves. Therefore, $\DDelta_0$-invariants coincide with the
invariants of $S$-equivalence. Note that the Seifert form of a knot
is a complete invariant of $S$-equivalence; in particular,
$S$-trivial knots are exactly the knots with trivial Alexander
polynomial.

One may ask which Vassiliev invariants are also
$\DDelta_0$-invariants. It was shown by Murakami and Ohtsuki in
\cite{mo} that all $\Q$-valued Vassiliev invariants that are
invariants of $S$-equivalence are polynomials in the coefficients of
the Conway polynomial. Here we extend their result to primitive
Vassiliev invariants with values in any abelian group $A$ with no $2$-torsion.

To this end, we define a logarithm-like function
$$\zlog\colon 1+z\cdot\Z[[z]]\to z\cdot\Z[[z]]$$
which takes multiplication to addition. This is, as far as we know, a novel construction which should be quite useful in studying primitive Vassiliev invariants
over the integers.
 Applying $\zlog$ to the Conway
polynomial $C(z)$, the coefficients are non-zero
only for even degrees of $z$. Denote by $pc_{2n}$ the coefficient at
$z^{2n}$.

\begin{thm}\label{dd0}
For each $n>0$,  $pc_{2n}$ is both a $\DDelta_0$-invariant and a
primitive Vassiliev invariant. Every primitive Vassiliev knot
invariant taking values in an Abelian group with no $2$-torsion, which is also a $\DDelta_0$-invariant, is a linear
combination of the $pc_{2n}$. As a Vassiliev invariant, $pc_{2n}$
has order $2n$, and $pc_{2n}\mod{2}$ has order $2n-1$.
\end{thm}

Here we should clarify what it means for an $A$-valued invariant $f$
to be a linear combination of the $pc_{2n}$: this is said to be the
case if there exists a homomorphism $\phi: \Z^{k}\to A$ for some
$k$, such that $f$ is the composition of the direct sum of the
$pc_{2n}$, for $0<n\leq k$, with $\phi$.

Proceeding to $\DDelta_1$ invariants, we show
\begin{thm}\label{dd1}

\noindent
\begin{itemize}
\item The associated graded $\Q$-module of
primitive rational-valued Vassiliev invariants which are also
$\DDelta_1$-invariants, has exactly one generator in degree $2n+1$
for each $n\geq 1$, and no generators in even degrees.

\item Any integer-valued Vassiliev invariant of order 4 is,
modulo $2$, a $\DDelta_1$-invariant.
\end{itemize}
\end{thm}

It is possible that there are more $2$-torsion invariants in larger
even degrees; we shall see, however that at most one copy of
$\Z/2\Z$ can exist in each even degree.

One could assemble the $\mathbb Z$-valued invariants into a power
series which is a $\DDelta_1$-invariant but not a Vassiliev
invariant, but we leave open the question of whether there exist
$\DDelta_1$-invariants which do not come directly from Vassiliev
invariants in this way.

We shall not discuss $\DDelta_n$-invariants with $n>1$ in detail,
though not for the lack of examples. Indeed, the $\DDelta$-move is a
special case of a null-move that preserves the ambient $3$-manifold
(see \cite{gk, gr}). Above we have mentioned that the rational lift
of the Kontsevich integral, $Z^{\mathfrak{rat}}$, is a universal
rational null-finite-type invariant of $S$-trivial knots. This
implies the following theorem.

\begin{thm}\label{zrat}
$Z_{i+1}^{\mathfrak{rat}}$ is a $\DDelta_{2i}$-invariant. Therefore,
the vector space of rational $\DDelta_{2n}$-invariants modulo
$\DDelta_{2n+1}$-invariants is infinite-dimensional.
\end{thm}

\subsection*{Some notation}
The band-pass and the doubled-delta moves will be denoted by BP and
$\DDelta$ respectively; $V$ will sometimes be used for the crossing
change. (Here ``$V$" stands for ``Vassiliev invariants".) The notation $\Z [X]$ will be used for the free
$\Z$-algebra on a monoid $X$. The lower central series of a group
$H$ will be denoted by $\gamma_kH$, indexed so that $\gamma_1H=H$. $\p{N}$ will denote the pure
braid group on $N$ strands; the trivial braid on $N$ strands will be
written as $\one_N$.

   %%%%%%%%%%%%%%%%%%%%%%%%%%%%%%%%%%%%%%%%%%%%%%%%%%%%%%%
   % dos.                                                %
   %                         %%%%                        %
   %                       %%    %%                      %
   %                      %%      %%                     %
   %                      %%      %%                     %
   %                              %%                     %
   %                             %%                      %
   %                            %%                       %
   %                          %%                         %
   %                        %%                           %
   %                      %%%%%%%%%%%                    %
   %                                                     %
   %%%%%%%%%%%%%%%%%%%%%%%%%%%%%%%%%%%%%%%%%%%%%%%%%%%%%%%

\section{Local Moves and Finite-Type invariants}\label{section:moves}
There are at least three different ways to define local moves on
knots. The definitions we give below are not equivalent; however,
all ``interesting" local moves, such as the band-pass, doubled
delta, $C_k$-moves and many other moves can be obtained from all
three constructions.

\begin{figure}
\begin{center}
\begin{tabular}{p{.40\linewidth}}
\begin{center}
\begin{minipage}{1.2in}\includegraphics[width=1.1in]{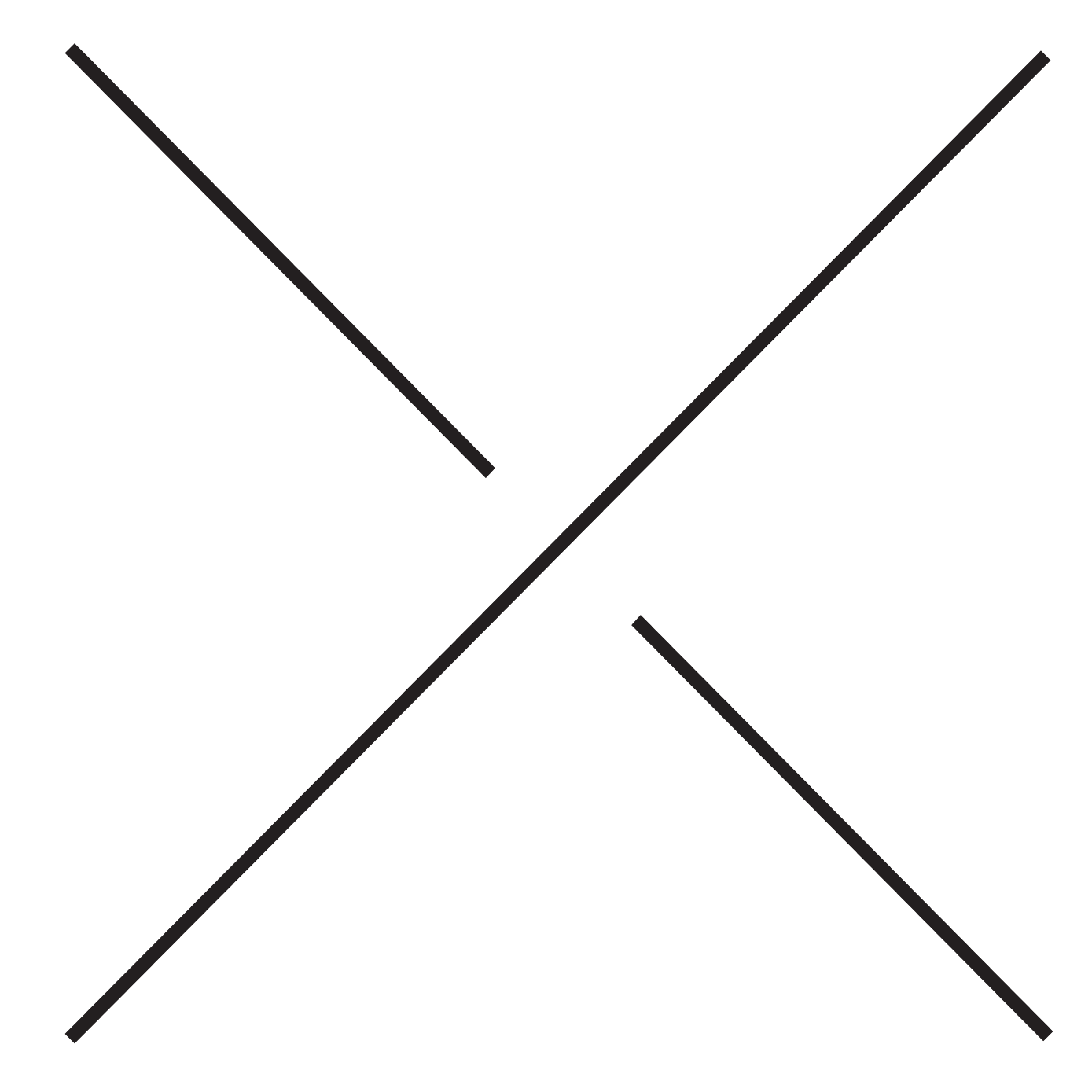}\end{minipage}$\leftrightarrow$ 
\begin{minipage}{1.2in}\includegraphics[width=1.1in]{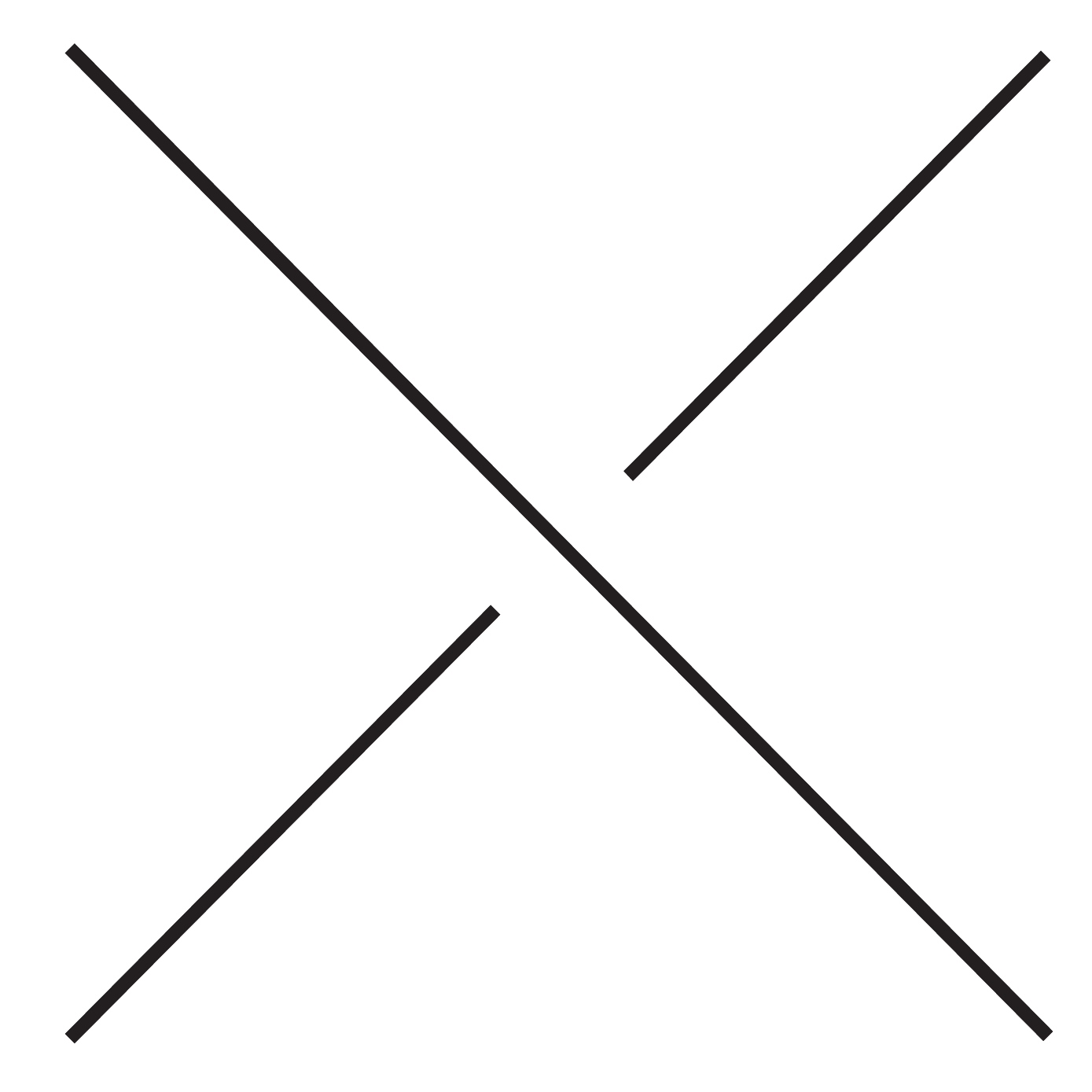}\end{minipage}\\
\end{center}\\
\\
\begin{center}Crossing change  $V$\end{center}
\end{tabular}\hfill
\begin{tabular}{p{.40\linewidth}}
\begin{center}
\begin{minipage}{1.2in}\includegraphics[width=1.1in]{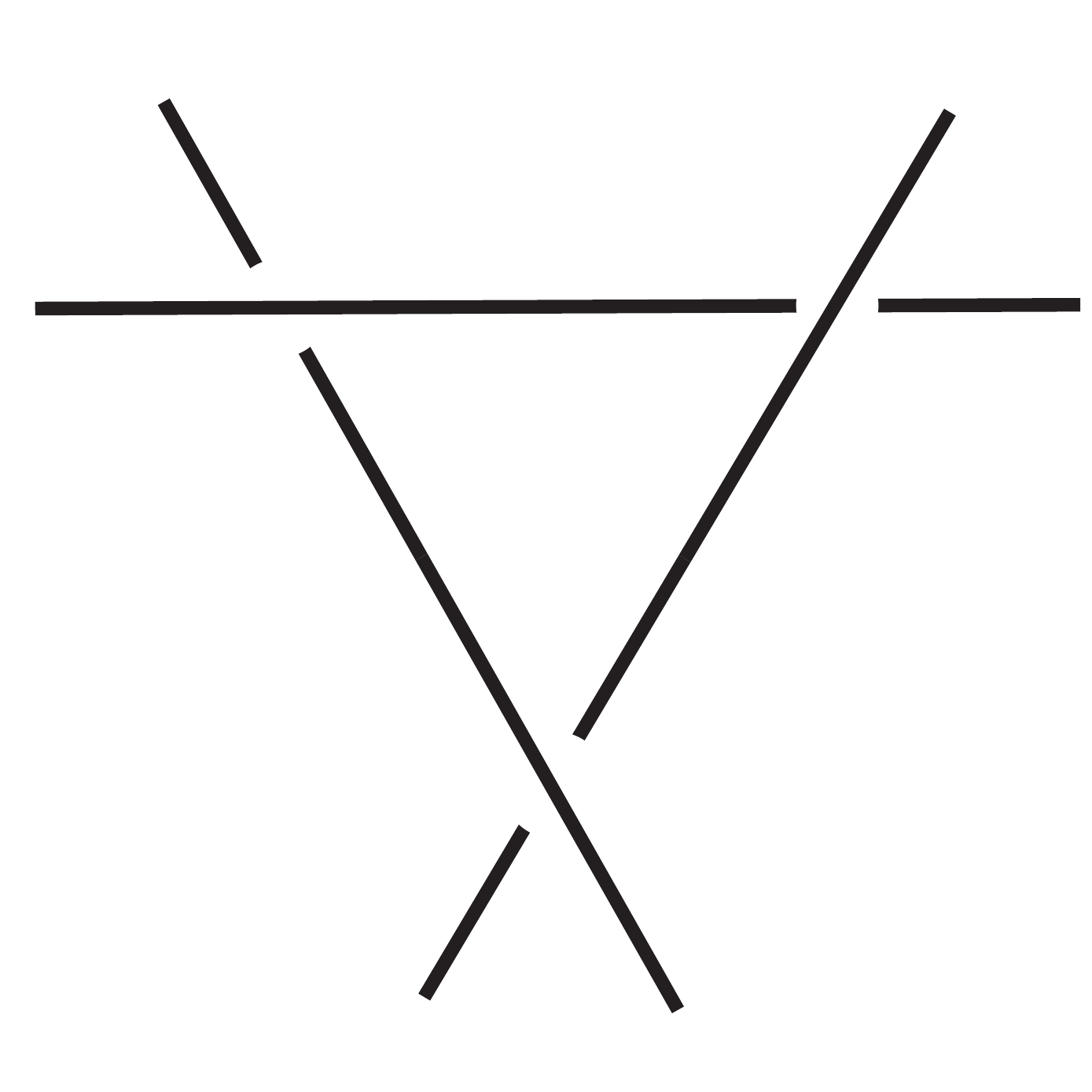}\end{minipage}$\leftrightarrow$ 
\begin{minipage}{1.2in}\includegraphics[width=1.1in]{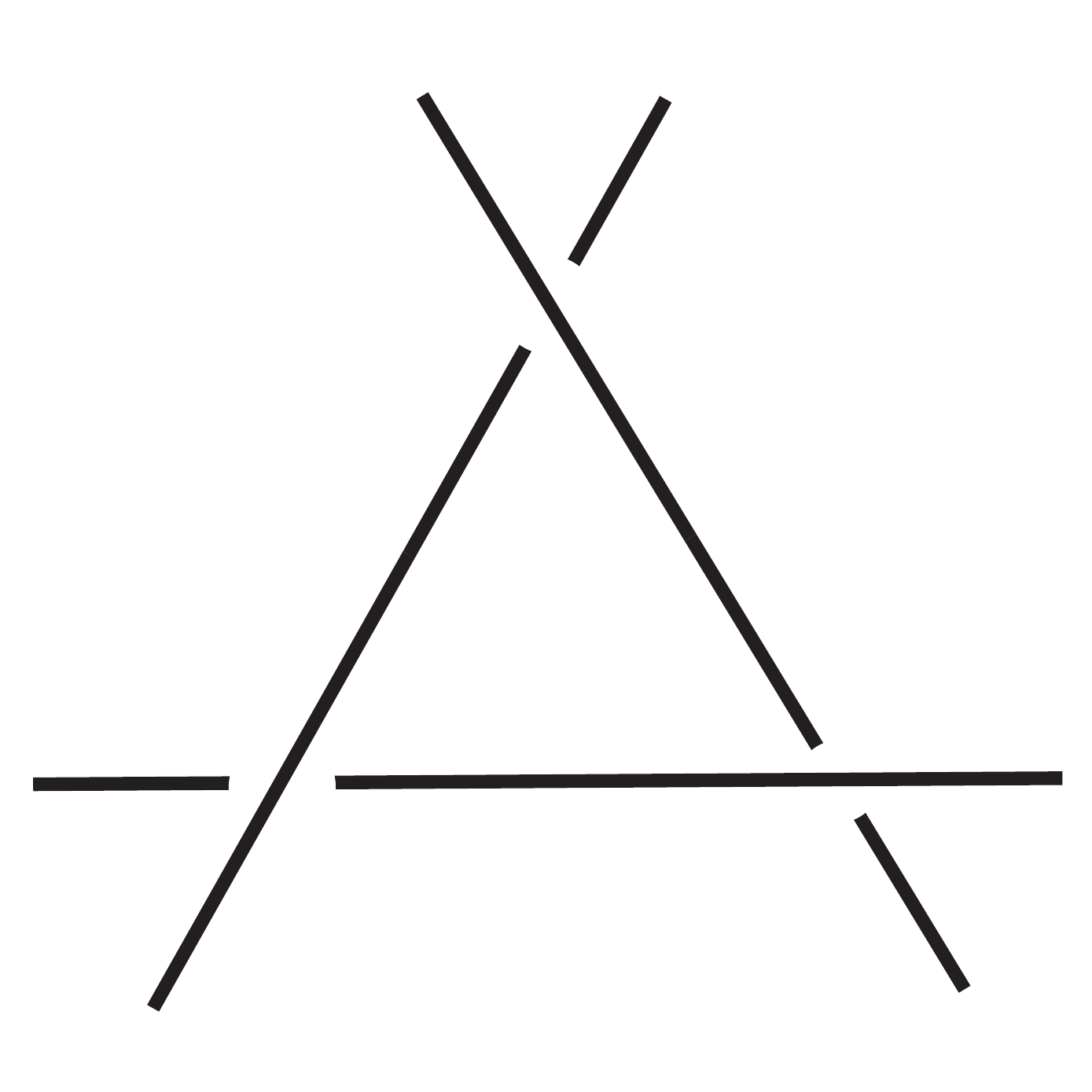}\end{minipage}
\end{center}
\\
\begin{center}Delta move  $\Delta$\end{center}
\end{tabular}\\
\begin{tabular}{p{.45\linewidth}}
\begin{center}
\begin{minipage}{1.2in}\includegraphics[width=1.1in]{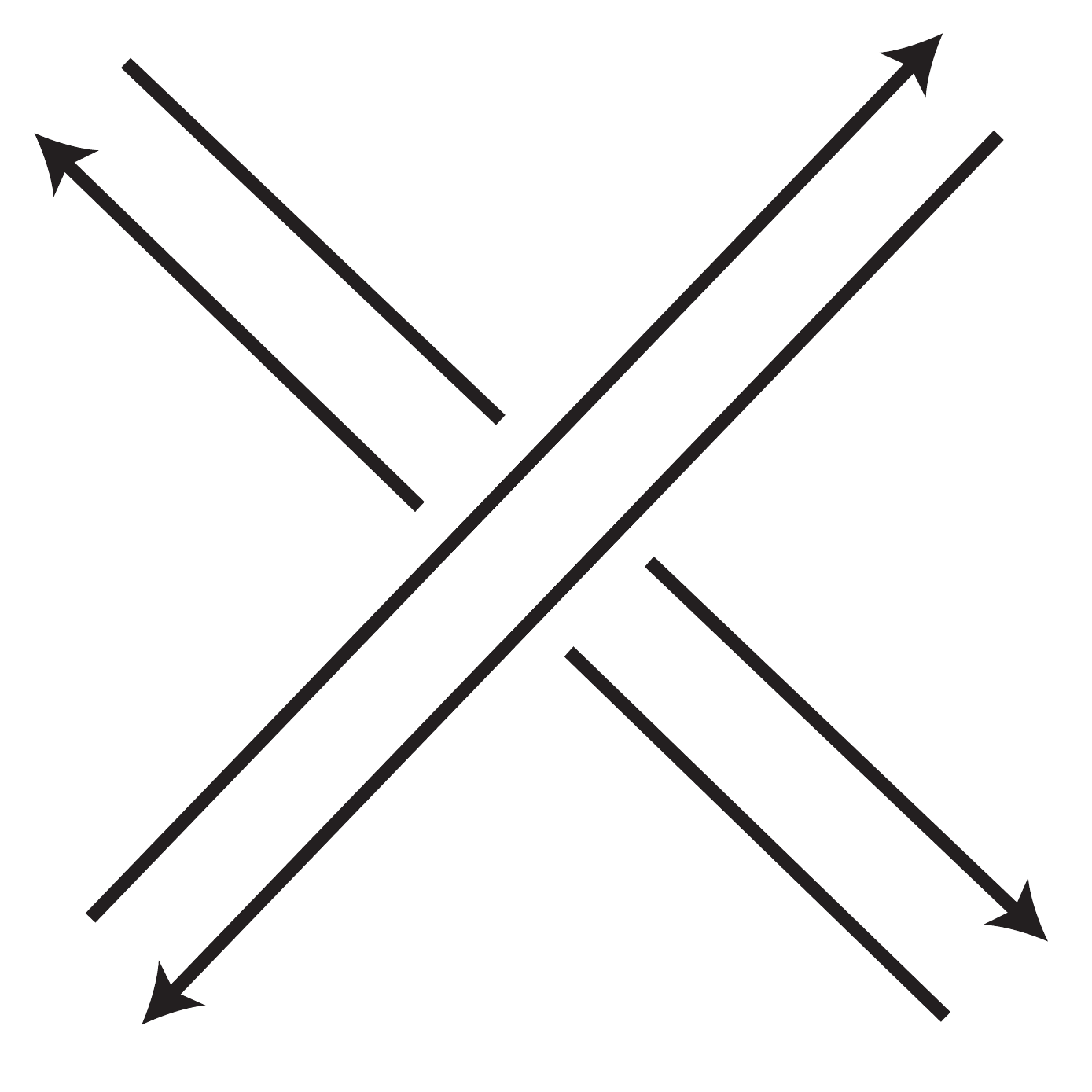}\end{minipage}$\leftrightarrow$ 
\begin{minipage}{1.2in}\includegraphics[width=1.1in]{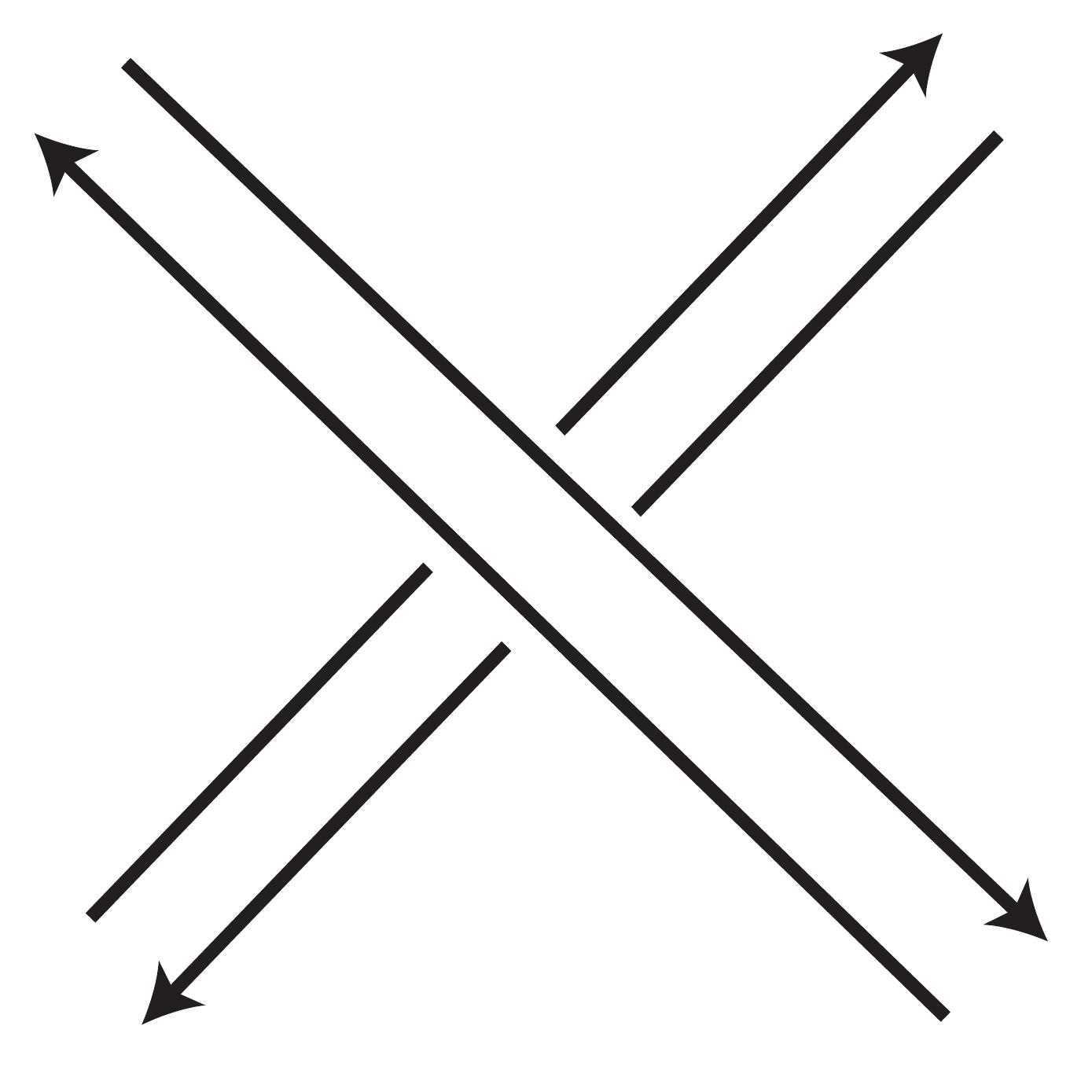}\end{minipage}\\
\end{center}\\
\begin{center}Band-pass move $\BP$ \end{center}
\end{tabular}\hfill
\begin{tabular}{p{.45\linewidth}}
\begin{center}
\begin{minipage}{1.2in}\includegraphics[width=1.1in]{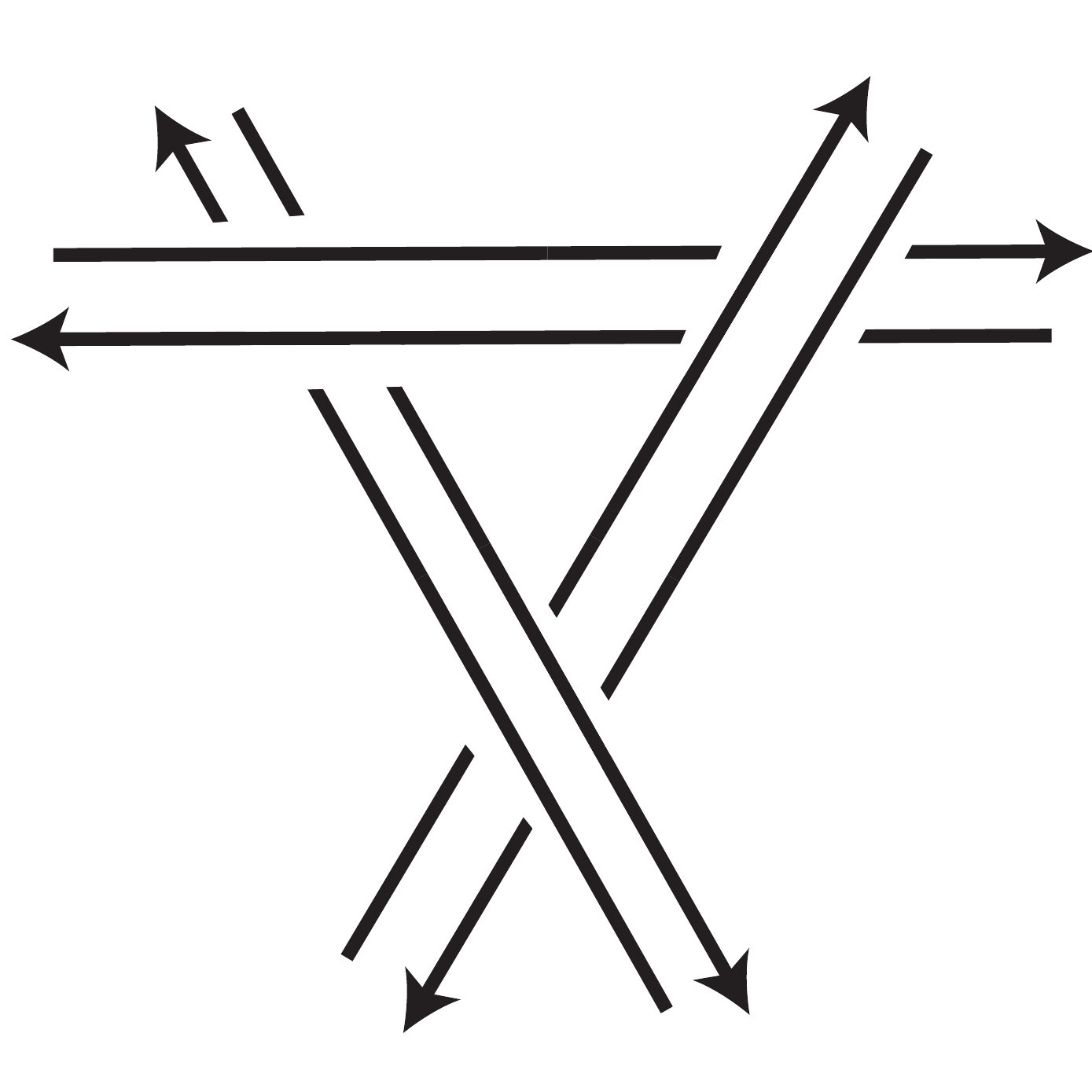}\end{minipage}$\leftrightarrow$ 
\begin{minipage}{1.2in}\includegraphics[width=1.1in]{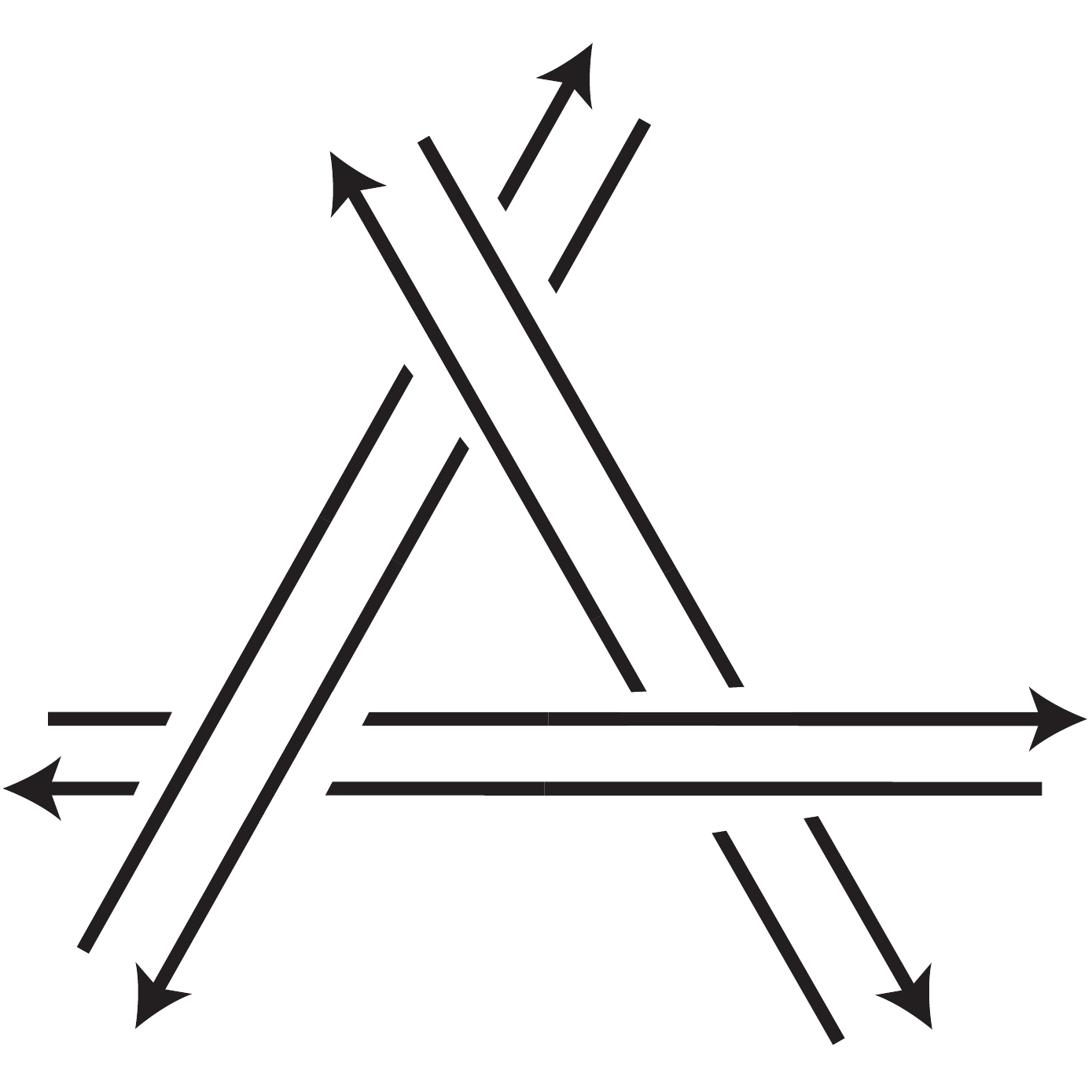}\end{minipage}\\
\end{center}\\
\begin{center}Doubled-delta move $\DDelta$ \end{center}
\end{tabular}\hfill
\\
\caption{Basic local moves on knots.}~\label{movetable}
\label{Moveset0}
\end{center}
\end{figure}

\subsection{Local moves and tangles.}
The most intuitive approach to local moves is via tangles. The
following definition can be found, for example, in \cite{ty2}.

Let $T$ and $S$ be two tangles in an oriented closed ball $U$.
Assume that neither $T$ nor $S$ contain closed components and denote
by $t_1, t_2, \ldots , t_k$ and $s_1, s_2, \ldots , s_k$ the arcs of
$T$ and $S$ respectively. Suppose that for each $t_i$ there exists
some $s_j$ such that $\partial t_i = \partial s_j$ . Then the
ordered pair $(T,S)$ is called  {\em a local move}. Two local moves
$(T,S)$ and $(T', S')$, defined in 3-balls $U$ and $U'$
respectively, are equivalent if there exists an orientation
preserving homeomorphism $h : U\to U'$ such that $h(T) = T'$ and
$h(S)$ is ambient isotopic to $S'$ relative to $\partial U'$.

Let $K_1$ and $K_2$ be knots in an oriented 3-manifold $M$. We say
that $K_2$ {\em is obtained from $K_1$ by applying a local move
$(T,S)$} if there exists a ball $U$ in $M$ such that
\begin{enumerate}
\item $\partial U$ intersects both knots transversely,
\item  $K_1$ and $K_2$ coincide outside $U$
\item the pair $(U\cap K_1, U\cap K_2)$ is a local
move equivalent to $(T,S)$.
\end{enumerate}

A theory of local moves based on this definition is developed in
\cite{ty1,ty2}. We shall not give further details here.

\subsection{Local moves and braids.}
\begin{figure}
$$
\xymatrix{
\xy
\xygraph{
     !{0;/r1.0pc/:} 
  !{\vcrossneg}
  [rru]!{\xcapv[1]@(0)}      
  [l]!{\vcross}    
  [lu]!{\xcapv[1]@(0)}  
  [r]!{\vcross}    
  [lu]!{\xcapv[1]@(0)}  
  !{\vcross}
  [rru]!{\xcapv[1]@(0)}    
  [ll]!{\vcross}
  [rru]!{\xcapv[1]@(0)}  
  [l]!{\vcrossneg}
  [lu]!{\xcapv[1]@(0)} 
  [r]!{\vcrossneg}
  [lu]!{\xcapv[1]@(0)} 
  !{\vcrossneg}
  [rru]!{\xcapv[1]@(0)}  
 [ruuuuuuuu] !{\xcapv[1]@(0)}
  !{\xcapv[1]@(0)}
 !{\xcapv[1]@(0)}
 !{\xcapv[1]@(0)}
 !{\xcapv[1]@(0)}
 !{\xcapv[1]@(0)}
 !{\xcapv[1]@(0)}
 !{\xcapv[1]@(0)}
 [ruuuuuuuu] !{\xcapv[1]@(0)}
 !{\xcapv[1]@(0)}
 !{\xcapv[1]@(0)}
 !{\xcapv[1]@(0)}
 !{\xcapv[1]@(0)}
 !{\xcapv[1]@(0)}
 !{\xcapv[1]@(0)}
 !{\xcapv[1]@(0)}
     }
     \endxy
     \ar@{|->}[r]
& 
\xy
\xygraph{
     !{0;/r1.0pc/:} 
  {\xcapv[2]@(0)<<}
  [ddr]!{\vcrossneg}
  [rru]!{\xcapv[1]@(0)}      
  [l]!{\vcross}    
  [lu]!{\xcapv[1]@(0)}  
  [r]!{\vcross}    
  [lu]!{\xcapv[1]@(0)}  
  !{\vcross}
  [rru]!{\xcapv[1]@(0)}    
  [ll]!{\vcross}
  [rru]!{\xcapv[1]@(0)}  
  [l]!{\vcrossneg}
  [lu]!{\xcapv[1]@(0)} 
  [r]!{\vcrossneg}
  [lu]!{\xcapv[1]@(0)} 
  !{\vcrossneg}
  [rru]!{\xcapv[1]@(0)}  
  [ruuuuuuuu] !{\xcapv[1]@(0)}
  !{\xcapv[1]@(0)}
 !{\xcapv[1]@(0)}
 !{\xcapv[1]@(0)}
 !{\xcapv[1]@(0)}
 !{\xcapv[1]@(0)}
 !{\xcapv[1]@(0)}
 !{\xcapv[1]@(0)}
 [ruuuuuuuu] !{\xcapv[1]@(0)}
 !{\xcapv[1]@(0)}
 !{\xcapv[1]@(0)}
 !{\xcapv[1]@(0)}
 !{\xcapv[1]@(0)}
 !{\xcapv[1]@(0)}
 !{\xcapv[1]@(0)}
 !{\xcapv[1]@(0)}
 [uuuuuuuulll]!{\vcap}
 [rr]!{\vcap}
 [ddddddddlll]!{\vcap-}
 [rr]!{\vcap-}
 [rr]!{\xcapv[2]@(0)<<}
     }\endxy
}
 $$
 \caption{The short-circuit closure of a braid.}~\label{shortcircuit}
\end{figure}
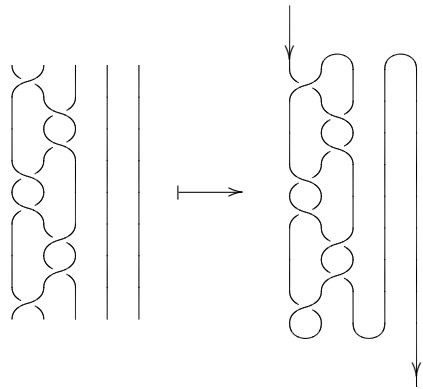

There are several kinds of braid closures. We shall use the
``short-circuit" closure of \cite{ms}; this is a version of the plat
closure which sends pure braids on an odd number of strands to
(long) knots. See Figure~\ref{shortcircuit}.
The short-circuit closure commutes with adding two
unbraided strands to a braid, hence, there is a map
$$\s :\p{\infty}\to \K.$$
A Markov-type theorem for the short-circuit closure is proved in
\cite{ms}: it says that there are two subgroups $H^T$ and $H^B$ of
$\p{\infty}$ such that the map $\s$ identifies the set of isotopy
classes of knots with the two-sided quotient
$H^T\backslash\p{\infty}/H^B$.

Let $\subgroup\subset\p{\infty}$ be a normal subgroup. Two knots $a$ and
$b$ are $\subgroup$-{\em equivalent} if there exist
$x,y\in\p{\infty}$ and $h\in\subgroup$ such that $a=\s{(x)}$,
$b=\s{(y)}$ and $x=hy$. The Markov Theorem for the
short-circuit closure implies that $\subgroup$-equivalence is indeed
an equivalence relation on the set of isotopy classes of knots. All
knots are $\p{\infty}$-equivalent.

An informal explanation of this definition in terms of local moves
is that $\subgroup$ should be thought of as the subgroup of
$\p{\infty}$, consisting of the braids that can be undone by some
fixed local move. For example, crossing changes can undo all braids
in $\p{\infty}$; delta moves undo the commutator subgroup
$[\p{\infty},\p{\infty}]$ (see \cite{st-delta}). Then one can think
of $\subgroup$-equivalent knots as those that can be transformed
into each other by a sequence of local moves that undo $\subgroup$.

Related to $\subgroup$-equivalence is the notion of a
$\subgroup$-finite-type invariant. Let $I_\subgroup$ be the
augmentation ideal of $\subgroup$ inside
$\Z[\subgroup]$, and denote by $\widehat{I}_{\subgroup}$ the ideal
in $\Z[\p{\infty}]$ generated by $I_{\subgroup}$. The ideal
$\widehat{I}_{\subgroup}$ is the kernel of the ring homomorphism
$\Z[\p{\infty}]\to\Z[\p{\infty}/\subgroup]$ sending elements of
$\p{\infty}$ to $\p{\infty}/\subgroup$.

A linear function $\Z[\K]\to A$, where $A$ is an abelian group, is a
{\em $\subgroup$-finite-type invariant of order $n$} if it vanishes
on all elements of $\Z[\K]$ which are of the form $\s{(a)}$ with %changed from $\s{ax}$
$a\in \widehat{I}_\subgroup^{n+1}$. For $\subgroup=\p{\infty}$ this
defines the usual Vassiliev invariants. Notice that all the knots
within the linear combination $\s{(a)}$ are $\subgroup$-equivalent,
so it makes sense to define a $\subgroup$-finite-type invariant on a
single $\subgroup$-equivalence class. In particular,
$\subgroup$-finite-type invariants of order $n$ for $\subgroup$-trivial knots can
be defined as those that vanish on $\s{(I_\subgroup^{n+1})}$.

Denote by $\tau_0:\p{\infty}\to\p{\infty}$ the homomorphism of
shifting the braid by two strands ``to the right". In other words,
$\tau_0$ sends the braid $x$ to $\one_2\otimes x$. For $k>0$ let
$\tau_k$ be the homomorphism of $\p{\infty}$ into itself that
triples the $k$-th strand.

\begin{thm}\label{thm:main}
Let $\m$ be a normal subgroup of $\p{\infty}$ such that
$\tau_k(\m)\subset \m$ for all $k\geqslant 0$. Then
\begin{enumerate}
\item $\m$-trivial knots considered modulo $\gamma_n\m$-equivalence form a
group under connected sum;
\item two $\m$-trivial knots are
$\gamma_{n}\m$-equivalent if and only if they cannot be
distinguished by $\m_{<n}$-invariants.
\end{enumerate}
\end{thm}

The proof of this theorem will be given in
Section~\ref{section:braids}.

\subsection{Local moves and claspers.}

We shall assume that the reader is familiar with the language of
claspers. For definitions and properties of claspers we refer to
\cite{habiro, ct}. We shall use the terminology of \cite{ct}.

Let $T$ be an abstract (i.e. not embedded) clasper of a fixed tree type.
Then a basic $T$-move on a pair $(\mfld,\knot)$ where $\knot$ is a
knot embedded in a three-manifold $\mfld$, is a surgery on an
embedded $T$-clasper in $\mfld\setminus \knot$, which produces a new
pair $(\mfld',\knot')$. A $cT$-move is a $T$-clasper surgery where
the leaves bound disjoint disks (called caps) which may hit the
knot. Note that such moves preserve the ambient $3$-manifold. An
$nT$-move is a $T$-move where the clasper's leaves link the knot
homologically trivially. (Here ``$n$'' stands for ``null.") Note that this makes
sense when $\mfld$ is a homology sphere. Finally an $ncT$-move is a
$T$-move where the leaves bound disjoint caps that intersect the
knot algebraically trivially.

Let the clasper with tree type an interval be denoted $I$, and the
clasper with tree type a ``wye" be denoted $Y$. Then what follows is
a table of other names for the above moves. See Figure~\ref{movetable}.

\begin{center}
\begin{tabular}{p{2em}|p{10em}p{10em}p{12em}}
&$c$&$n$&$nc$\\
\hline
\\
$I$&Crossing change (V)&\ &Band-pass move (BP)\\

$Y$&Delta move ($\Delta$) &Null-move ($\Null$)&Doubled-delta move ($\DDelta$)\\

\end{tabular}
\end{center}

\medskip

This table should be interpreted in the sense that every $cI$-move
can be realized by a \emph{sequence} of crossing changes and
similarly for the other entries.

Let $\move$ denote a move type of the form $cT$ or $ncT$. Then there
is  a descending filtration
$$
\Z[\K]=\F^\move_0\supset\cdots\supset
\F^\move_n\supset \F^\move_{n+1}\supset\cdots,
$$
where $\F_n^\move$ is defined as follows. Given a knot $\knot$, and
$n$ disjointly embedded $\move$-claspers $C_1,\ldots, C_n$, with
disjointly embedded caps. Let
$$
[\knot;C_1,\ldots, C_n]=\sum_{\sigma\subset
[n]}(-1)^{|\sigma|}\knot_{\sigma},
$$
where $\knot_\sigma$ means ``modify $\knot$ by the claspers $C_i$"
where $i\in\sigma$. Then $\F_n^\move$ is defined to be the submodule
generated by such elements $[\knot;C_1,\ldots, C_n]$ for all choices
of $\knot$ and $\{C_i\}$. By Habiro's zip construction  and standard arguments, in the $cT$ case
the submodule $\F_n^\move$ is actually generated by such alternating sums
where each clasper is simple. (That is, each of the caps intersect the knot
in a single point.) In the $ncT$ case one may assume each cap intersects the knot in a pair
of algebraically canceling points.

Let $\overline{\Z[\K]}$ be the quotient of $\Z[\K]$ by the subspace
generated by the vectors $$\knot_1+ \knot_2-\knot_1\# \knot_2$$ for
all pairs of knots. Let $\overline{\F^\move_n}\subset
\overline{\Z[\K]}$ be the image of ${\F^\move_n}$. Linear functions
$\overline{\Z[\K]}\to G$ are exactly the primitive (that is,
additive) knot invariants.

\begin{defi}\label{def:ft-claspers}

\noindent\begin{enumerate}
\item An invariant $f\colon \K\to A$, where $A$ is an abelian group,
is said to be {\em $\move$-finite-type of degree (order) $n$} if
$f(\F^\move_{n+1})=0$. We also say that $f$ is a {\em
$\move_n$-invariant}.
\item Two knots $\knot_1$ and $\knot_2$ are {\em $\move_n$-equivalent}
if $\knot_1-\knot_2\in \F^\move_{n+1}$. In other words,
$\move_n$-equivalent knots are those that cannot be distinguished by
$\move_{\leq n}$ invariants. By $\mugroup_n$ we denote the set of
knots which are $\move_0$-equivalent to the unknot, considered
modulo $\move_n$-equivalence.
\end{enumerate}
\end{defi}

Theorem~\ref{thm:main} implies the following statement.
\begin{thm}\label{thm:group}
Given a tree $T$ which is a $Y$ or an $I$, let $\move$ be a move of the type $cT$ or $ncT$.
Then $\mugroup_n$ is an abelian group for all $n$.
\end{thm}

The language of claspers facilitates the treatment of the null-move
$\Null$. One can define a null-finite-type invariant of degree $n$
to be an invariant (with values in some abelian group) of pairs
$(\mfld, \knot)$ which vanishes on alternating sums of $n+1$ null
moves, and analogously define $\Null_n$-equivalence. Let
$G^{\Null}_n$ be the set of pairs $(\mfld,\knot)$ which are
$\Null_0$-equivalent to the pair $(\sphere,\unknot)$, considered
modulo $\Null_n$-equivalence. Note that $G^{\Null}_n$ is a monoid
under the operation of connected sum, where the $3$-ball that is
removed hits the knot in a standard arc. Let
${G}^{\Null}_n(\sphere)$ be the submonoid of pairs
$(\sphere,\knot)$.

\begin{thm}
${G}^{\Null}_n(\sphere)$ is a group for all $n$.
\end{thm}

\begin{proof}
Let $X$ denote the set of pairs $(S^3,\knot)$ which are
$\Null_0$-equivalent to the pair $(\sphere,\unknot)$. We claim that
$X$ is the set of all pairs $(\sphere, \knot)$ where $\knot$ has
trivial Alexander polynomial. Every such pair is in $X$ by
\cite{ns}. Conversely, any element in $X$ will have trivial
Alexander polynomial.

But then $X$ is precisely the set of knots in $\sphere$ which are
$\DDelta_0$ equivalent to the unknot. $G^{\DDelta}_n$ is an abelian
group, by Theorem~\ref{thm:group}, which projects onto the monoid
${G}^{\Null}_n(\sphere)$. The latter is then forced to be a group.
\end{proof}

\begin{problem}
Is $G^{\Null}_n$ a group?
\end{problem}

On the level of manifolds, this is true. The monoid of homology
$3$-spheres is a group after modding out by alternating sums of $Y$
moves as above \cite{gouss2,habiro}.

\subsection{Translating between claspers and braids.}

Given a tree $T$, a {\em simple $cT$ clasper} (with respect to a
given knot) is defined to be a $cT$ clasper where each leaf bounds a
disk which hits the knot in one point. 

\begin{defi}
A clasper $C$ embedded in the complement of the trivial braid $\one_n$ is said to be \emph{braid-like}, if
\begin{enumerate}
\item Surgery along the clasper $(\one_n)_C$ is a braid.
\item Regarding $\one_n$ as a subset of a cube in the usual way, the complement of $\one_n$ in the cube is a regular neighborhood of the clasper.
\end{enumerate}
\end{defi}

\begin{lem}
If $T$ is an interval or a $Y$, then there is a braid-like $T$ clasper. 
\end{lem}
\begin{proof}
The case of an interval is obvious. A braid-like $Y$-clasper is illustrated below.
\begin{center}
\begin{minipage}{1.4in}
\includegraphics[width=1.3in]{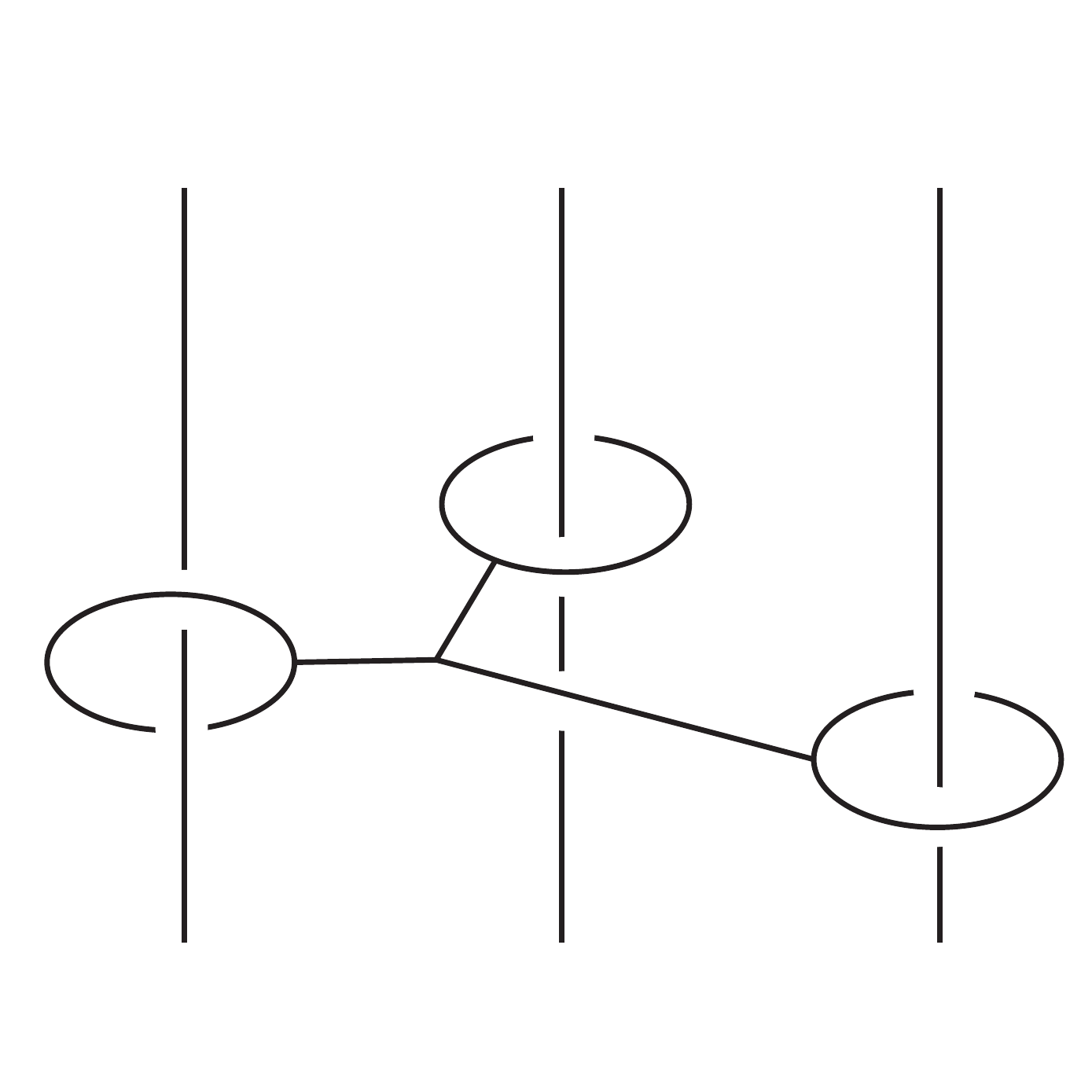}
\end{minipage}
 $\leftrightarrow $
 \begin{minipage}{1.4in}
\includegraphics[width=1.3in]{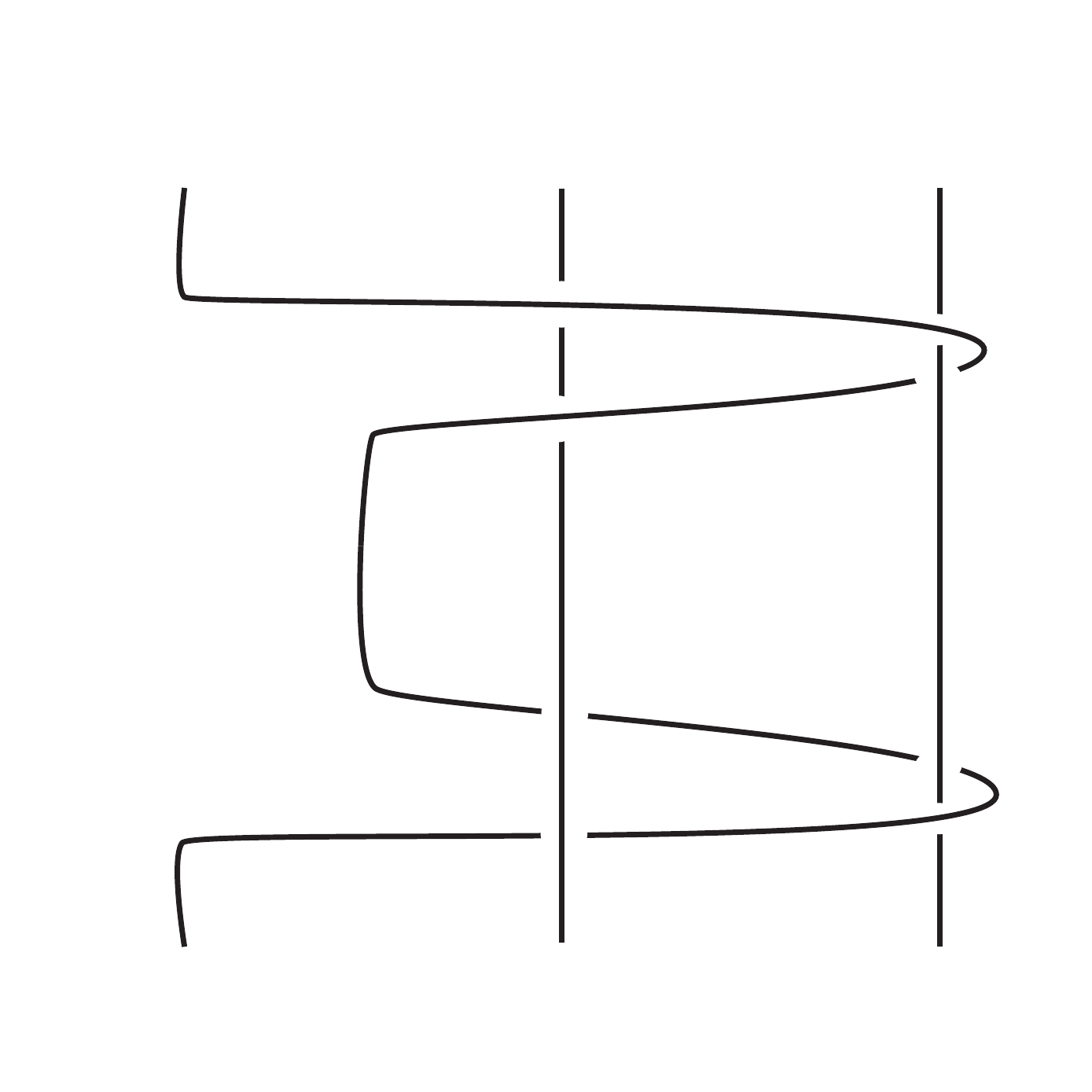}
\end{minipage}
\end{center}

\end{proof}

\begin{defi}\item
\begin{enumerate}
\item Let $\subgroup^{cT}$ be the subgroup of $\p{\infty}$ normally
generated by elements $(\one_N)_{C_1\cup\cdots\cup C_i}$, where $C_1\cup\cdots\cup C_i$ is a union of
$cT$ claspers which have the property that $(\one_N)_{C_1\cup\cdots \cup C_i}$ is a braid. 
\item Denoting by $\double\colon 
\p{\infty}\to\p{\infty}$ the homomorphism which doubles every
strand, define $\subgroup^{ncT}$ as the group normally generated
by $\double(\subgroup^{cT})$ inside $B_\infty$, the direct limit of
braid groups on a finite number of strands. Note that
$\subgroup^{ncT}$ is indeed a subgroup of $\p{\infty}$. 
\end{enumerate}
\end{defi}

\begin{prop}
The 
 subgroups $\subgroup^{ncT}$ and $\subgroup^{cT}$ are invariant under the
 strand tripling homomorphism $\tau_k$.
\end{prop}
\begin{proof}
Consider a generating element $(\one_N)_{C_1\cup\cdots\cup C_i}$ of $\subgroup^{cT}$. If we triple the $k$th strand, this will either miss the caps of the claspers or it will convert a single intersection point to three intersection points. The claspers will remain capped. Also, the homomorphism $\tau_k$  commutes with the clasper surgery, indicating that  $\tau_k((\one_N)_{C_1\cup\cdots\cup C_i})$ is again a braid. We note that the inverse of a generator is also a generator, since it is formed by mirror reflection, and the mirror reflection of a clasper surgery is also a clasper surgery of the same tree type. (Note though that the clasper itself is not the mirror reflection, since orientation data is important.)

Now consider a generator $\double((\one_N)_{C_1\cup\cdots\cup C_i})$ of $\subgroup^{ncT}.$  Applying $\tau_k$, we are tripling one strand of the double of a strand in $\double((\one_N)_{C_1\cup\cdots\cup C_i})$. That is we are doubling all the strands of the braid except one which we are quadrupling. So    $\tau_k\double((\one_N)_{C_1\cup\cdots\cup C_i})=
\double((\one_{N+1})_{C'_1\cup\cdots\cup C'_i})$ where the claspers $C'_i$ are formed from $C_i$ by doubling the appropriate strand of $\one_N$.
\end{proof}

\begin{thm}\label{jim}
Suppose $T$ is a tree for which a braid-like clasper exists.
For $\mu=cT$ or $\mu=ncT$, $\move_n$-invariants coincide with
$\subgroup^\move_n$-invariants.
\end{thm}
\begin{proof}
Let us first consider the case $\mu=cT$. Suppose $f\colon\K\to G$ is
a $cT_n$-invariant. Note that $\hat{I}_{\subgroup^{cT}}^{n+1}$ is
generated by elements of the form $$g_1(h_1-1)g_2(h_2-1)\cdots
g_{n+1}(h_{n+1}-1)g_{n+2},$$ where $h_i\in \subgroup^{cT}$ and
$g_i\in\p{\infty}$. One may assume that each $h_i$ is a normal
generator of $\subgroup^{cT}$: $h_i=\one_{\overline{C}_i}$, where $\overline{C}_i$ is a union
of $cT$ claspers.

Letting
$$x=\s\left( g_1(\one_{\overline{C}_1}-1)g_2(\one_{\overline{C}_2}-1)\cdots
g_{n+1}(\one_{\overline{C}_{n+1}}-1)g_{n+2}\right)$$
we must show that
$f(x)=0.$ But this follows since
$$x=(-1)^{n+1}[\s(g_1\ldots g_{n+2}); \overline{C}_1,\ldots,\overline{C}_{n+1}],$$
and in the usual way one can write an alternating sum of surgeries on unions of claspers
as a linear combination of alternating sums of surgeries on individual claspers.

Conversely, let $C_1,\ldots, C_{n+1}$ be a union of simple $cT$
claspers on a knot $\knot$.  We want to deform the Morse function on
$\mathbb R^3$ to get a knot $\s(b)$ where the claspers $C_i$ satisfy the following properties:
\begin{itemize}
\item Each $C_i$ sits in an interval $[s_i,t_i]$ with respect to the Morse function, and all these intervals are pairwise disjoint.
\item The restriction of the knot to each interval is the trivial braid.
\item Each $C_i$ is braid-like. 
\end{itemize}
To do this, using condition (2) in the definition of ``braid-like," position the claspers $C_i$ 
so that they match their braid-like embedding with respect to the morse function for $\one_n$.
Now poke the maxima and minima of the knot up and down, to get a short-circuit closure. 
  Since a
surgery on a braid-like clasper is equivalent to a modification by
an element of $\subgroup^{cT}$, the case $\mu=cT$ is settled.

Consider now the case $\mu=ncT$. Showing that an $ncT_n$ invariant vanishes on $\s(\hat{I}^{n+1}_{\subgroup^{ncT}})$ is proved exactly analogously. We have that $\hat{I}^{n+1}_{\subgroup^{ncT}}$ is  generated
by elements
$$g_1(h_1-1)g_2(h_2-1)\cdots g_{n+1}(h_{n+1}-1)g_{n+2}$$
where $h_i\in \double(\subgroup^{cT})$, and $g_i\in B_\infty$. So
each $h_i$ is equal to $\one_{\overline{C}_i}$, where $\overline{C}_i$ is a union of
$ncT$-claspers, and the argument indeed proceeds as above.

Now for the converse. Let $C_1,\ldots, C_{n+1}$ be a union of $ncT$
claspers on $\knot$ where each leaf bounds a disk hitting $\knot$ in
precisely two algebraically canceling points. As above, we can
deform the Morse function so that each $C_i$ sits in between two heights $s_i,t_i$, 
all the intervals $[s_i,t_i]$ are disjoint,  in each interval
 $C_i$ sits in the complement of a trivial braid such that $C_i$
 is obtained from a braid-like embedding by doubling some strands.
 So each $C_i$ converts the identity  element $\one_N$ to a
pure braid $b$. Assume that the strands of $\one_N$ are numbered
from left to right. The leaves of $C_i$ each grab two antiparallel
strands. If this pair of strands is of the form $(2j-1,2j)$ we have
$\one_{C_i}=\double(\one_{C^\prime_i})$ with $C'_i$ a braid-like
clasper. If this is not the case, let $p$ be a braid such that in $p
\one_{C_i}p^{-1}$, the leaves of $C_i$ do grab pairs of consecutive
strands of the form $(2j-1,2j)$. So then $C_i$ turns $\one_N=p^{-1}p
\one_N p^{-1}p$ into $b'_i=p^{-1}(pbp^{-1})p$. Now
$pbp^{-1}=\double(\one_{C^\prime_i})$ for a braid-like $cT$-clasper
$C^\prime_i$, and so $b'_i\in \subgroup^{ncT}$. As before, the
argument is now straightforward. Assume, for notational ease, that
the claspers $C_i$ occur in order on the braid which closes up to
give $\knot$. Write
$$\knot=\s(g_1\cdot \one\cdot g_2\cdots \one\cdot g_{n+2})$$
where $C_i$ turns the $\one$ following $g_i$ into the braid
$b'_i\in\subgroup^{ncT}$. Then $$[\knot;C_1,\ldots,
C_{n+1}]=\s(g_1(\one-b'_1)\cdots g_{n+1}(\one-b'_{n+1})g_{n+2})\in
\hat{I}_{\subgroup^{ncT}}^{n+1}.$$
\end{proof}

%Sometimes the subgroups of $\p{\infty}$ corresponding to a given
%clasper can be explicitly identified. We have already
%mentioned that if $T$ is the $Y$-graph, the corresponding $cT$-move
% the delta move whose corresponding subgroup of $\p{\infty}$ is
%just the commutator subgroup $\gamma_2\p{\infty}$. 

%[SKETCH OF ARGUMENT:
%$\mathcal G^{cY}=\gamma_2\p{\infty}$. $\subset$ is fairly easy to prove. 
%To see $\supset$, consider a commutator in the braid group. Then the first strand cobounds
%a capped surface with the trivial first strand. This capped surface yields a union of capped 
%$Y$ claspers, by the zip construction, we may assume they are a union of simple $Y$ claspers. 
%Should we allow arbitrary claspers in defn of $\mathcal G^{cY}$?
%]

   %%%%%%%%%%%%%%%%%%%%%%%%%%%%%%%%%%%%%%%%%%%%%%%%%%%%%%%
   % tres.                                               %
   %                         %%%%                        %
   %                       %%    %%                      %
   %                      %%      %%                     %
   %                              %%                     %
   %                             %%                      %
   %                           %%                        %
   %                             %%                      %
   %                      %%      %%                     %
   %                      %%      %%                     %
   %                       %%    %%                      %
   %                         %%%%                        %
   %                                                     %
   %%%%%%%%%%%%%%%%%%%%%%%%%%%%%%%%%%%%%%%%%%%%%%%%%%%%%%%

\section{Groups of knots via the short-circuit closure}
\label{section:braids}

\subsection{Proof of Theorem~\ref{thm:main}.}

We shall think of $\p{\infty}$ as consisting of braids on finite odd
number of strands. Given $x\in\p{\infty}$, the long knot $\s(x)$ is
obtained by joining together lower and upper ends of the strands in
turns, and extending to infinity the upper end of the first strand
and the lower end of the last strand. There are two subgroups $H^T,
H^B\subset \p{\infty}$ such that  $x,y\in\p{\infty}$ satisfy
$\s(x)=\s(y)$ if and only if there exist $t\in H^T$ and $b\in H^B$
with $x=tyb$, see \cite{ms}. %We shall not need the explicit
%definition of the subgroups $H^T$ and $H^B$ here.

Recall the definition of the endomorphisms $\tau_i$ of $\p{\infty}$:
$\tau_0$ shifts the braids by two strands ``to the right" and, for
$k>0$, $\tau_k$ triples the $k$-th strand. It is clear that
$\s{(\tau_0{(x)})}=\s{(x)}$ for all $x\in\p{\infty}$, thus
$\tau_0{(x)}=txb$, where $t\in H^T$ and $b\in H^B$.
\begin{lem}\label{important}
For $x\in \m$ there exist $t\in H^T\cap \m$ and $b\in H^B\cap \m$
such that $\tau_0{(x)}=txb$.
\end{lem}
\begin{proof}
Let $t_{2k-1}=\tau_{2k-1}(x)(\tau_{2k}(x))^{-1}$, and let
$b_{2k}=(\tau_{2k+1}(x))^{-1}\tau_{2k}(x)$. Notice that $t_{2k-1}\in
H^T\cap \m$ and $b_{2k}\in H^B\cap \m$. We have
$$ \tau_{2k-1}(x)=t_{2k-1}\tau_{2k}(x),$$
$$\tau_{2k}(x)=\tau_{2k+1}(x)b_{2k}.$$
There exists $N$ such that $\tau_{2N+1}(x)=x$. Thus the following
equality holds:
$$\tau_0(x)=t_{1}\cdots t_{2N-1}xb_{2N}\cdots b_{0},$$
and this completes the proof.
\end{proof}

Denote by $G_n$ the set of $\gamma_{n}\m$-equivalence classes of
$\m$-trivial knots. The operation of  connected sum descends to an
addition on $G_n$ as $\tau_0(\gamma_{n}\m)\subset\gamma_{n}\m$.
\begin{lem}
$G_n$ is an abelian group.
\end{lem}
\begin{proof}
The connected sum of knots is commutative and associative and the
neutral element is provided by the equivalence class of the unknot.
Thus it is only necessary to check that every element of $G_n$ has
an inverse.

Let $\knot=\s{(y)}$ be the short-circuit closure of
$y\in\gamma_{m}(\m\cap\p{2N-1})$. If $m\geqslant n$, the knot
$\knot$ is $\gamma_{n}$-equivalent to the trivial knot and defines
the identity in $G_n$. Assume that $m<n$. Consider the knot
$\knot\#\s{(y^{-1})}=\s{(y\tau_0^{N}(y^{-1}))}$. By
Lemma~\ref{important} the braid $\tau_0^{N}(y^{-1}))$ can be written
as $ty^{-1}b$ with $t\in H^T\cap \m$ and $b\in H^b\cap \m$, so
\[ y\tau_0^{N}(y^{-1}) = yty^{-1}b = t(t^{-1}yty^{-1})b \]
and
\[\s{(y\tau_0^{N}y^{-1})}=\s(t(t^{-1}yty^{-1})b)=
\s(t^{-1}yty^{-1})=\s{(y_{1})} \] where $y_{1}\in
\gamma_{m+1}\m\cap\p{4N-1}$.

Repeating the construction  one sees that
$\knot\#\s{(y^{-1})}\#\s{(y_{1}^{-1})}\#\ldots\#\s{(y_{q-1}^{-1})}$
is $\gamma_{(m+q)}\m$-equivalent to the unknot. In particular, the
knot $\s{(y^{-1})}\#\ldots\#\s{(y_{n-m-1}^{-1})}$ is an inverse for
$\knot$ in $G_n$.
\end{proof}
In fact, the same proof works to show that $\m$-trivial knots modulo
the {\em derived} series of $\m$ form a group.

\begin{problem}
Are these groups modulo the derived series related to known knot
invariants? Do they contain any concordance information?
\end{problem}

\medskip

The value of any $\m_{<n}$-invariant on a knot depends only on the
$\gamma_{n}\m$-equivalence class of the knot. Indeed,
$\gamma_{n}\m-\one\subset I^n_\m$.

In order to prove Theorem~\ref{thm:main} it remains to show that
different $\gamma_{n}\mathcal G$-equivalence classes can be distinguished by
finite-type invariants. This is established by the following lemma.

\begin{lem}\label{lemma:last}
The quotient map $\s{(\m)}\to G_{k}$ is a $\m_k$-invariant.
\end{lem}

\begin{proof}

We have to show that $\s{(I_\m^{n})}$ is mapped to zero in
$G_{n-1}$.

Define a {\em relator of degree $n$ and length $m$} as an element of
$\Z[\s(\m)]$ of the form
$$\s((x_1-1)(x_2-1)\ldots(x_m-1)y)\leqno{(\ast)}$$
with $y\in \m$, $x_i\in\gamma_{n_i}\m$ and $\sum n_i = n$. The
greatest $n$ such that a relator is of degree $n$ will be called the
{\em exact degree} of a relator. A {\em composite relator} is an
element of $\Z[\s(\m)]$ of the form
$\knot_1\#\knot_2-\knot_1-\knot_2+1$ where $\knot_1,\knot_2$ are
knots. Notice that a connected sum of two relators of non-zero
degree is a linear combination of composite relators.

The kernel of the map $\Z[\s(\m)]\to G_n$ contains all the relators
of length 1 and degree $n$ and all the composite relators. On the
other hand, an element of $\s{(I_\m^n)}$ is a linear combination of
relators of length $n$ and, hence, of degree $n$. Thus if we show
that any relator of degree $n$ is a linear combination of relators
of degree $n$ and length $1$ and composite relators, the lemma is
proved.

Suppose that there exist relators of degree $n$ which cannot be
represented as linear combinations of the above form. 
Among such
relators choose the relator $R$ of minimal length and, given the
length, of maximal exact degree. ( If one of the $n_i$ in a relator exceeds $n$,
then clearly the relator is in the kernel of the quotient map, so there is indeed an upper bound on the exact degree
of such a relator.)

Assume that $R$ is of the form $(\ast)$ as above, with $y,
x_i\in\p{2N-1}\cap \m$. Choose $t\in H^T$ and $b\in H^B$ such that
the braid $tx_{1}b$ coincides with the braid obtained from $x_{1}$
by shifting it  by $2N$ strands to the right, that is, with
$\tau_0^N(x_1)$. By Lemma~\ref{important} the braids $t$ and $b$ can
be taken to belong to the same term of the lower central series of
$\m$ as the braid $x_1$. The relator
\[ R'=\s((tx_{1}b-1)(x_{2}-1)\ldots(x_{m}-1)y) \]
is a connected sum of two relators and, hence, is a combination of
composite relators. On the other hand,
\[\begin{array}{lcl}
R'-R & = & \s((tx_{1}b-x_{1})(x_{2}-1)\ldots(x_{m}-1)y)\\
     & = & \s(x_{1}(b-1)(x_{2}-1)\ldots(x_m-1)y)
\end{array}\]
Notice now that $(b-1)$ can be exchanged with $(x_i-1)$ and $y$
modulo relators of shorter length or higher degree. Indeed,
\[(b-1)y=y(b-1)+([b,y]-1)yb\]
and
\[(b-1)(x_i-1)=(x_i-1)(b-1)+ ([b,x_i]-1)(x_{i}b-1)+([b,x_{i}]-1).\]
Thus, modulo  relators of shorter length or higher degree
\[\s(x_{1}(b-1)(x_{2}-1)\ldots(x_{m}-1)y)=
\s(x_{1}(x_{2}-1)\ldots(x_{m}-1)y(b-1))=0.\] and this means that $R$
is a linear combination of composite relators and relators of length
1 and degree $n$.
\end{proof}

   %%%%%%%%%%%%%%%%%%%%%%%%%%%%%%%%%%%%%%%%%%%%%%%%%%%%%%%
   % cuatro.                                             %
   %                            %%%                      %
   %                           %%%%                      %
   %                          %% %%                      %
   %                         %%  %%                      %
   %                        %%   %%                      %
   %                       %%    %%                      %
   %                      %%     %%                      %
   %                     %%      %%                      %
   %                    %%%%%%%%%%%%%                    %
   %                             %%                      %
   %                             %%                      %
   %                                                     %
   %%%%%%%%%%%%%%%%%%%%%%%%%%%%%%%%%%%%%%%%%%%%%%%%%%%%%%%

\section{Band-pass invariants}
In this section we analyze band-pass invariants, completely
characterizing the primitive ones.

\begin{thm}\label{bpthm}
For $n\geq 1$, primitive $\BP_{n}$-invariants  coincide with
primitive Vassiliev invariants of order $n$. The unique primitive
$\BP_{0}$-invariant is the Arf invariant. In other words,
$G^{\BP}_n$ is the index $2$ subgroup of $G^{V}_n$ consisting of
those knots with trivial Arf invariant.
\end{thm}

We first prove some lemmas.

\begin{lem}\label{a1}

\noindent\begin{enumerate}
\item \label{exper} For all $n$, $\F^{\BP}_{n}\subset \V{n}$
\item \label{e2} For $n>2$, $\overline{\V{n}}\subset \overline{\F^{\BP}_{n-2}}$.
\end{enumerate}
\end{lem}
\begin{proof}
Part (\ref{exper}) follows by writing a band-pass move as a union of
crossing changes.

For part(\ref{e2}), we use \cite[Prop. 6.10]{habiro}. The quoted
result implies that $\F^{V}_{n}$ is generated by elements of the
form $\knot-\knot_C$ where $C$ is a simple (which implies capped)
tree clasper of degree $n$ on a knot $\knot$, and by elements which
lie in the second power of the augmentation ideal. However, these
latter elements vanish upon passing to $\overline{\Z[\K]}$.
 By
\cite[Thm. 28]{ct}, we can assume the underlying tree of $C$ is a
tree where every trivalent vertex shares an edge with at least one
univalent vertex.

Consider such a tree clasper $C$ on a knot $\knot$. Decompose $C$ in
the following way: insert two Hopf pairs of leaves along all edges
which are not incident to a univalent vertex. This gives us a union
of $I$ and $Y$ claspers. Perform a surgery of $\knot$ along the $Y$
claspers to get a knot $\knot_Y$. Note that $\knot_Y$ is isotopic to
$\knot$. There are exactly $n-2$  of the $I$ claspers, denote them by
$I_1,\ldots, I_{n-2}$.  Then
$$(-1)^{n-2}[\knot_Y;I_1,\ldots,I_{n-2}]=\knot_C-\knot.$$
Since the $I_j$ are of the form $ncI$, $\knot_C-\knot\in
\F^{\BP}_{n-2},$ completing the proof.
\end{proof}

\begin{lem}\label{a2}

\noindent\begin{enumerate}
\item \label{a21} For $n>2$, $\overline{\F^{\BP}_n}$ surjects
onto $\overline{\V{n}}/\overline{\V{n+1}}$.
\item \label{a22} $\overline{\F^{\BP}_2}$ surjects onto
$2\Z\subset \Z\cong\overline{\V{2}}/\overline{\V{3}}$.
\end{enumerate}
\end{lem}
\begin{proof}
Assume first $n>2$. Then $ \V{n}/\V{n+1}$ is generated by chord
diagrams with $n$ chords (see, for example, \cite{bn}). By
\cite{ng}, this is generated by wheels attached to the Wilson loop
by some permutation, and by separated diagrams. Passing to $
\overline{\V{n}}/\overline{\V{n+1}}$, the separated diagrams vanish.
Thus it suffices to realize a wheel with $n$ legs attached to the
Wilson loop by some alternating sum of band-pass moves. An easy way
to do this is to realize the wheel by an embedded clasper, producing
band-pass moves as in the proof of Lemma~\ref{a1}~(\ref{e2}).

If $n=2$ the same argument shows that the subspace generated by the
wheel with two legs is realized. This is indeed equal to twice the
generator. Finally, the generator itself cannot be realized since it
has nontrivial Arf invariant, and Arf is a band-pass degree zero
invariant \cite{k}, and so vanishes on $\overline{\F^{\BP}_2}$.
\end{proof}

\begin{prop}\label{a3}
For $n>2$ we have $\overline{\F^{\BP}_n}=\overline{\V{n}}$.
\end{prop}
\begin{proof}
By Lemma~\ref{a1}~(\ref{exper}), it suffices to prove the inclusion
$\overline{\V{n}}\subset\overline{\F^{\BP}_n}$. Let
$v_n\in\overline{\V{n}}$. By Lemma~\ref{a2}~(\ref{a21}), there is
some $v_{n+1}\in\overline{\V{n+1}}$ such that $v_n+v_{n+1}\in
\overline{\F^{\BP}_n}$. Again by Lemma~\ref{a2}, there is some
$v_{n+2}\in\overline{\V{n+2}}$ such that $v_{n+1}+v_{n+2}\in
\overline{\F^{\BP}_{n+1}}\subset\overline{\F^{\BP}_n}$. But by
Lemma~\ref{a2}~(\ref{e2}) $v_{n+2}\in \overline{\F^{\BP}_n}$.
Therefore $v_n\in\overline{\F^{\BP}_n}$ as desired.
\end{proof}

\begin{proof}[Proof of Theorem~\ref{bpthm}]
The first statement is a direct corollary of Proposition~\ref{a3}.
The fact that the unique degree zero invariant is the Arf invariant
was proven by Kauffman \cite{k}. Now suppose $f$ is a primitive
$\BP$-finite-type invariant of degree $1$. Then it is of Vassiliev
degree $3$. Now
$\overline{\F^{\BP}_2}\supset\overline{\F^{\BP}_3}\twoheadrightarrow\overline{\V{3}}/\overline{\V{4}}$,
implying that $f$ is actually of Vassiliev degree $2$. But
$\overline{\F^{\BP}_2}\twoheadrightarrow 2\Z\subset
\Z\cong\overline{\V{2}}/\overline{\V{3}}$. Thus $f$ is either zero
or the mod 2 reduction of the degree $2$ Vassiliev invariant. That
is, it must be zero or the Arf invariant.
\end{proof}

By the theory of Hopf algebras, the set of rational-valued finite-type invariants is a polynomial algebra over the primitive invariants. The same may not be true for $\mathbb Z$-valued invariants, leaving open the possibility that there are non-primitive $\BP$-finite-type invariants which are not finite-type in the standard sense. Thus we pose the following problem.
\begin{problem}
Are $\BP$-finite-type invariants always products of primitive
invariants?
\end{problem}

   %%%%%%%%%%%%%%%%%%%%%%%%%%%%%%%%%%%%%%%%%%%%%%%%%%%%%%%
   % cinco.                                              %
   %                      %%%%%%%%%%                     %
   %                      %%      %%                     %
   %                      %%                             %
   %                      %% %%%%                        %
   %                      %%     %%                      %
   %                              %%                     %
   %                              %%                     %
   %                      %%      %%                     %
   %                      %%      %%                     %
   %                       %%    %%                      %
   %                         %%%%                        %
   %                                                     %
   %%%%%%%%%%%%%%%%%%%%%%%%%%%%%%%%%%%%%%%%%%%%%%%%%%%%%%%

\section{doubled-delta invariants}
In this section, we explore $\DDelta_n$-invariants.
$\DDelta_0$-equivalence is precisely the same as $S$-equivalence,
the condition that two knots have Seifert surfaces with identical
Seifert forms, see \cite{ns}. However, the question still arises as
to which Vassiliev invariants are also $\DDelta_0$-invariants.
Murakami and Ohtsuki \cite{mo} have shown that in the case of
$\Q$-valued invariants, they must be polynomials in $c_{2n}$, the
coefficients of the Conway polynomial. We prove the analogous result
for a more general target group, namely abelian groups with no $2$-torsion, but restricting to primitive
invariants. (Corollary \ref{Sequivintro}). This still leaves open the
following question:

\begin{problem}
Do all (not necessarily primitive) Vassiliev invariants which are
also $S$-equivalence invariants come from the Conway polynomial?
(This is true for $\Q$-valued invariants by \cite{mo}, and for
primitive invariants taking values in an abelian group with no $2$-torsion by Corollary \ref{Sequivintro}.)
\end{problem}

Moving on to the analysis of primitive $\DDelta_1$-invariants, we
show that for each $n\geq 1$ there exists exactly one integer-valued
Vassiliev invariant of order $2n+1$ which is also a
$\DDelta_1$-invariant, and that there are no integer-valued
$\DDelta_1$-invariants that are Vassiliev of even order. On the
other hand, we shall see that any integer-valued Vassiliev invariant
of order 4 is, modulo $2$, a $\DDelta_1$-invariant.

There may be more $2$-torsion invariants in larger even degrees, but
there are two stumbling blocks to proving this. First, we are not
sure whether the 2-torsion exists on the level of weight systems,
although we can show that no more than one copy of $\Z/2\Z$ exists
in each even degree. The second stumbling block is more serious.
Once we have a weight system, we do not have a mod $2$ Kontsevich
integral to produce a knot invariant, which is how we prove the
first part of Theorem~\ref{dd1}.

Finally, as observed in the introduction, the Euler degree $n+1$
part of the Kontsevich integral, $Z^{\mathfrak{rat}}_{n+1}$, is a
surjective $\DDelta_{2n}$-invariant. Thus the kernels
$$\ker\left(\group{\DDelta}{2n}\twoheadrightarrow\group{\DDelta}{2n-1}\right)$$
are infinitely generated for all $n\geq 1$.

\subsection{$\DDelta_0$-equivalence and diagrams}

Let $\A_n$ be the $\Z$-module of {\em diagrams} (called ``Chinese
character diagrams" in \cite{bn}). These are connected trivalent
graphs with a distinguished oriented cycle ({\em Wilson loop}, or
{\em outer circle}) and a cyclic ordering of edges at each vertex
not on the Wilson loop. The edges that are not part of the Wilson
loop form the {\em internal graph} of the diagram.

Let $\A^I_n$ be $\A_n$ modulo separated diagrams \cite{ct2}. Any
primitive $\gv{n}$-invariant induces a weight system on $\A^I_n$.

An {\em insulated vertex} in a closed diagram of $\A_n$ is a vertex
which neither lies on the Wilson loop, nor is connected to it by an
edge. Let $\A^{iv}_n$ be the submodule of $\A^I_n$ generated by
those diagrams which have an insulated vertex.

\begin{prop}\label{iv}
Suppose $f$ is a Vassiliev invariant which is also a
$\DDelta_0$-invariant. Then the corresponding weight system
$\hat{f}$ vanishes on diagrams with insulated vertices.
\end{prop}

\begin{proof}
One way to prove this is through claspers. (Compare \cite{mo}.)
Given a diagram, $D$, with an insulated vertex, let $C_1,\ldots,
C_k$ be simple claspers on the unknot, which are embedded versions
of the connected components of the inner graph of the diagram.
%
%   A FIGURE MIGHT HELP HERE    :)
%
Then $\hat{f}(D)=\sum_{\sigma\subset \{ C_i\}}(-1)^{|\sigma|}
f(\unknot_\sigma)$. (See \cite[Lemma 2.1]{ct2}.) Let $C_{i_0}$ be
the clasper corresponding to the component with an insulated vertex.
Let $\knot$ be $\unknot$ modified by a fixed collection of $C_i$
except $C_{i_0}$. Then it suffices to show that surgery on the
component $C_{i_0}$ does not change the $S$-equivalence class of
$\knot$, because then all of the terms in the alternating sum cancel
in pairs.

A vertex, $v$, of a diagram $D$ is called {\em good} if $D\setminus
v$ is a connected graph (where we include the Wilson loop.)
Otherwise the vertex is called {\em bad}. Our strategy will be to
show that a diagram that has an insulated vertex is either zero or
has a good insulated vertex.

Suppose $D$ has a bad vertex. That means that there is a separating
edge in the diagram. Consider a separating edge which separates the
diagram into two pieces $G$ and $H$, where $H$ is the piece
containing the Wilson loop, and where $G$ does not itself contain a
separating edge. If $G$ has a trivalent vertex, then it is a good
insulated vertex. Otherwise $G$ is a loop, and so the diagram is
zero, since it has a loop edge\footnote{A diagram with a loop edge
is zero even over $\Z$.} .

Now, let $v_0$ be a good insulated vertex. Break the clasper
$C_{i_0}$ into a union of claspers by inserting Hopf pairs of
leaves on every edge incident to $v_0$. Let the $Y$ clasper which
contains $v_0$ be called $C_Y$, and let the union of the other
claspers be called $C_O$. Note that each of the claspers in
$C_O$ has a leaf on the knot. Thus
$\knot_{C_{i_0}}=(\knot_{C_O})_{C_Y}$. Now $C_Y$ is a capped
$Y$-clasper on the knot $\knot_{C_O}$ whose leaves link the knot
trivially. Therefore $C_Y$ can be realized as a union of doubled
delta moves. (To see this, use Habiro's zip construction to write the $Y$
as a union of capped $Y$ claspers each of which hits the knot in
two algebraically canceling points.)

The argument is completed by noticing that $\knot_{C_O}=\knot$ since
the claspers in $C_O$ each have a leaf bounding a disk disjoint from
the knot.
\end{proof}

Next, we shall construct a weight system which vanishes on diagrams
with insulated vertices. One way to construct weight systems is
through the intersection graph of a chord diagram. This graph is
defined by drawing a vertex for every chord, and putting an edge
between vertices if the chords cross. (Two chords are said to cross
if the four endpoints of the two chords alternate between the two
chords as you travel around the Wilson loop.)

A {\em Hamiltonian cycle} in the intersection graph is a cycle which
goes through every vertex exactly once without repeating any edges.

\begin{defi}
Let $\ws\colon \A_n\to \Z/2$ be the function defined as the number
of unoriented Hamiltonian cycles in the intersection graph, modulo $2$.
\end{defi}

\begin{lem}
The function $\ws$ is a primitive weight system.
\end{lem}
\begin{proof}
We need to show that $\ws$ vanishes on 4T relators and separated
diagrams, for $n\geq 4$.

It is clear that $\ws$ vanishes on separated diagrams since the
intersection graph of a separated diagram is not connected.

As for the $4T$ relation, let $\alpha$ and $\beta$ be the two chords
affected by our 4T relation. There are three complementary regions
of the Wilson loop. Let chords joining each of the three pairs be
called $A,B,C$. There are also chords that don't join two regions.
Call those $D$.  This is illustrated below:

\centerline {\includegraphics[width=.8\linewidth]{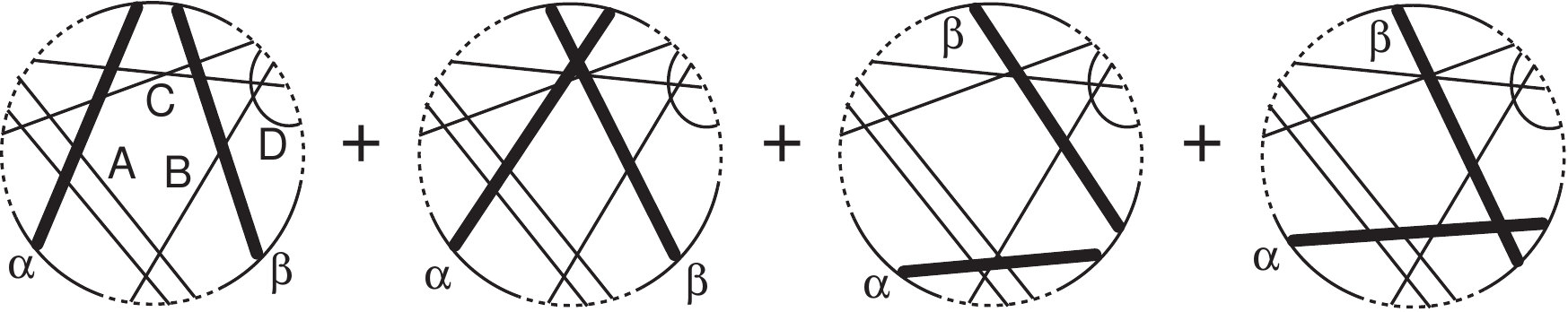} }

In the four terms of our relation we have the following
intersections between $\alpha$, $\beta$ and the other chords:
\begin{itemize}
\item[1:]$(\alpha, A),(\alpha, C),(\beta,B),(\beta,C)$

\item[2:] $(\alpha, A),(\alpha, C),(\beta,B),(\beta,C),
(\alpha,\beta)$

\item[3:]$(\alpha, A),(\alpha, B),(\beta,B),(\beta,C)$

\item[4:]$(\alpha, A),(\alpha, B),(\beta,B),(\beta,C),
(\alpha,\beta)$
\end{itemize}

We show that there are an even number of Hamiltonian cycles in these
four terms. First, there is evidently a $1-1$ correspondence between
cycles in the first intersection graph and cycles that don't go
through $(\alpha,\beta)$ in the second. (Similarly for third and
fourth.) Thus the total number of Hamiltonian cycles that don't go
through the bridge $(\alpha,\beta)$ is even. Next, we show there is
a $1-1$ correspondence between cycles that go through the bridge in
two and those that go through the bridge in four, {and that change
letters as they do so}. The first type of cycle must have a subchain
of the form $(A\text{ or }C,\alpha,\beta,B\text{ or }C)$ and the
second of the form $(A\text{ or }B,\alpha,\beta,B\text{ or }C)$. One
can easily check that the  input-output unordered pairs of different
letters are in correspondence in both cases. Finally, we show that
there is an even number of cycles that cross the bridge back into
the same letter, say $C$.  $\alpha$ and $\beta$ must be connected to
distinct vertices $c_1,c_2$ of $C$, else we would only have a
$3$-cycle. That means that the endpoints of the bridge each connect
to $c_1,c_2.$ If $c_1,\alpha,\beta, c_2$ is part of a Hamiltonian
cycle, so is $c_2,\alpha,\beta,c_1$, and we have shown that the
number of Hamiltonian cycles is even.
\end{proof}

\begin{lem}
$\ws=1$ on every wheel attached to the Wilson loop by some
permutation.
\end{lem}
\begin{proof}
Let $F_1$ and $F_2$ be two graphs with the same set of vertices, and
assume that $F_2$ is obtained from $F_1$ by adding one edge. Then we
can depict the formal sum $F_2+F_1$ by a graph that has the same
edges and vertices as $F_2$, but where the edge which is missing in
$F_1$ is dashed. Similarly, a graph with $n$ dashed edges is defined
as a sum of $2^n$ graphs.

Now, apply the STU relation to every leg of the wheel. The result is
a sum of $2^n$ chord diagrams (since we are working over $\Z/2\Z$
the signs are irrelevant); in all these diagrams each chord comes
from an edge of the circle part of the wheel. The corresponding
$2^n$ intersection graphs all share a common part, and there are $n$
edges, corresponding to adjacent edges in the wheel, which may or
may not be contained in the intersection graph. Hence, the sum of
the $2^n$ intersection graphs can be written as a graph with $n$
dashed edges. Note that the dashed edges in this graph form a
Hamiltonian cycle, and all other Hamiltonian cycles in this graph
come in pairs.

Indeed, if a Hamiltonian cycle does not go through some dashed edge,
then there are an even number of summands with this cycle, since the
edge being dashed means that we can either insert it or not. Thus we
are reduced to counting the Hamiltonian cycles that go through every
dashed edge, of which there is exactly one!
\end{proof}

\begin{lem}\label{onevanish}
$\ws$ vanishes on diagrams with at least one insulated vertex.
\end{lem}
\begin{proof}
Any diagram with an insulated vertex can be rewritten as a linear combination of
diagrams with an insulated vertex and connected internal graph, and
of separated diagrams. Since $\ws$ vanishes on the latter, we now
consider non-separated diagrams.

If a diagram has one loop and one insulated vertex, then it can be
reduced to an even number of wheels, by applying IHX to eliminate
trees growing off the loop.

If a diagram has two loops, then breaking one of the loops apart via
IHX and STU we can write it as a sum of an even number of diagrams with one
loop. These are either all wheels, or they are diagrams that have
one loop and also an insulated vertex.

In the former case, the fact that $\ws=1$ on wheels  finishes the
argument, whereas the latter case has already been covered.

Finally, if the diagram is a tree with an insulated vertex, then
modulo diagrams with one loop and one insulated vertex we can slide
the endpoints of the tree into any position. In particular we can
make two edges emanating from a trivalent vertex terminate in
adjacent positions on the Wilson loop. It is known (see, for example, \cite{ng})
that the resulting triangle formed by these two edges and the arc of
the Wilson loop joining them, can be collapsed to a point and
turned into a bubble inserted at any $3$-valent vertex.
(In other words, the trivalent vertex is replaced by a circle to which three edges attach.)
This puts us
back in the case of a diagram with a loop and an insulated vertex.
\end{proof}

\begin{prop}\label{calc}
If n is even $\A^I_n/\A_n^{iv}=\Z$, and if $n$ is odd then
$\A^I_n/\A_n^{iv}=\Z/2$.
\end{prop}
\begin{proof}
First, we look at $n\geq 3$. We know that $\A^I_n$ is generated by
wheels attached to the Wilson loop.
Modulo diagrams with two loops (which must have an insulated
vertex) the order in which the legs hit the Wilson loop is
irrelevant. If $n$ is odd, there is an orientation reversing
symmetry, so that $\A^I_n/\A_n^{iv}$ is at most $\Z/2$. However, the
weight system $\ws$ detects the wheel, so that we get $\Z/2$ on the
nose. When $n$ is even, the $\Z$ is detected by an integer lift of
the weight system $\ws$ \cite{bg}.

If $n=2$ the only diagrams with insulated vertices are already zero
modulo IHX, and the resulting space is just $\Z$.
\end{proof}

\subsection{$\DDelta_0$-invariants}
Propositions~\ref{iv} and \ref{calc} imply that for each degree there is at most
one primitive Vassiliev invariant which is also a $\DDelta_0$-invariant, since
such an invariant must give rise to a specific weight system.

In this section, we will show that there are primitive invariants
coming from the Conway polynomial that do indeed give rise to these
weight systems.

\subsubsection{A discrete version of the logarithm}
Define a function
$$\exp_{\Z}\colon x\cdot\Z[[x]]\to 1+x\cdot \Z[[x]]$$ by the rule
$\exp_{\Z}(\sum_{i}a_i x^i)=\prod_{i}(1+(-x)^i)^{a_i}$. Indeed $$\exp_{\Z}(\sum_{i}a_i x^i)=1-a_1x+\left(a_2+\binom{a_1}{2}\right)x^2-\left(a_3+a_1a_2+\binom{a_1}{3}\right)x^3+\cdots$$
In general, the coefficient of $x^\ell$ is given by $(-1)^\ell\displaystyle\sum_{i_1a_1+\cdots i_ka_k=\ell}\prod_{j=1}^k\binom{a_j}{i_j}$. In particular the coefficient of $x^\ell$ is equal to $\pm a_\ell$ plus a polynomial in the lower coefficients. This implies that $\exp_{\mathbb Z}$ is a bijection. 
Let $\zlog$ be the inverse function. 
\begin{lem}
Writing $\log_{\Z}(1+\sum_{i=1}^\infty b_ix^i)=\sum_{i=1}^\infty a_ix^i$,
 we have $a_i=(-1)^ib_i+q_i(b_1\ldots,b_{i-1})$ where $q_i$ is a polynomial of graded degree $i$.
\end{lem}
\begin{proof} 
We have
$$1+\sum_{i=1}^\infty b_ix^i=\exp_{\Z}(\sum_{i=1}^\infty a_ix^i)=1+\sum_{i=1}^\infty(-1)^i(a_i+Q_j(a_1,\ldots,a_{i-1}))x^i$$ where $Q_i$ is a polynomial of graded degree $i$. The equation $b_i=(-1)^j(a_i+Q_j(a_1,\ldots,a_{i-1}))$ can be recursively solved. Assume that all $a_j$ for $j<i$ are of the form $a_j=(-1)^jb_j+q_j(b_1,\ldots,b_{j-1})=(-1)^jb_j+q_j$. Then
$$a_i=(-1)^ib_i-(-1)^iQ_i\left(-b_1,\ldots,(-1)^{i-1}b_{i-1}+q_{i-1}\right)=(-1)^ib_i+q_i,$$
\end{proof}

Now define the invariants $pc_{2i}$ by the equation
$$\zlog C(z) =\sum_i pc_{2i} z^{2i}$$
where $z^2=x$ in the definition of $\zlog$. Recall that $C(z)$ is the Alexander-Conway polynomial.

Then
$$pc_{2i}=(-1)^ic_{2i}+q(c_2,c_4,\ldots,c_{2i-2}),$$ where $q$ is a
polynomial of graded degree $2i$. 
Thus these invariants $pc_{2i}$ are primitive Vassiliev invariants
of order $2i$. In fact, we have

\begin{prop}\label{nosystem}

\noindent\begin{enumerate}
\item $pc_{2n}$ is a  Vassiliev invariant of order $2n$ for all $n$.
\item $pc_{2n}\mod 2$ is a  Vassiliev invariant of order $2n-1$ for all $n\geq 2$.
\end{enumerate}
Moreover $\frac{1}{2}pc_{2n}$ induces the weight system
$\hat{\ws}\colon \A^I_{2n}\to \Z$ and $pc_{2n}\mod 2$ induces the
weight system $\ws\colon\A^I_{2n-1}\to \Z_2$. Here $\hat{\ws}$ is the integer lift of $\ws$ defined in \cite{bg}.
\end{prop}

\begin{cor}\label{Sequivintro}
Suppose $f\colon{\emph{\K}}\to A$ is a primitive Vassiliev invariant of order $k$ which is also a $\DDelta_0$-invariant. Suppose further that $A$ is an abelian group with
no $2$-torsion. Then $k$ is even and there is a
homomorphism $\phi\colon \Z^{\oplus k/2}\to A$ such that
$$f=\frac{1}{2^{k/2}}\phi\left(\bigoplus_{i=1}^{k/2} pc_{2i}\right).$$
\end{cor}

\begin{proof}[Proof of Corollary~\ref{Sequivintro}]
We prove this by induction on the  Vassiliev degree of $f$. Suppose
$f\colon\K\to A$ is a Vassiliev invariant of order $k$ which is primitive and invariant under the doubled-delta move. Let
$\hat{f}\colon \A^I_{k}\to A$ be the induced weight system. By Proposition~\ref{calc}, $k$ must be even
and there is a homomorphism
$\phi\colon\Z\to A$ such that $2\hat{f}= \phi\hat{\ws}$. So
$2f-\phi(pc_{k})$ is a Vassiliev invariant of order $k-1$ (and therefore of order $k-2$) which is
also a $\DDelta_0$ invariant. In this way,
we can write
$$2^{k/2}f=\phi_2(pc_2)+\phi_4(pc_4)+\cdots+ \phi_k(pc_k),$$
where $\phi_i$ is a homomorphism from $\Z$ to $A$. 
Tensoring with $\Z[\frac{1}{2}]$, which is an injection on $A$, we get  that $$f=\frac{1}{2^{k/2}}(\phi_2(pc_2)+\phi_4(pc_4)+\cdots +\phi_k(pc_k))\in A.$$
\end{proof}

This theorem leaves open the possibility that there are invariants taking values in an abelian $2$ group which do not come from the Conway polynomial.

\begin{problem}
Are there any $2$-torsion Vassiliev invariants which are also $\DDelta_0$ invariants, that do not arise from the Conway polynomial? 
\end{problem}

Before we prove Proposition~\ref{nosystem}, we need a lemma.

\begin{lem}\label{nosystem2}
The degree $2n$ weight system for $pc_{2n}$ is equal to $\pm 2$ on
wheels attached to the Wilson loop. The weight system for
${pc}_{2n}$ in degree $2n-1$ is odd on wheels. 
\end{lem}
\begin{proof}
First we note that the permutation by which the wheel is attached to
the Wilson loop is irrelevant. The difference of two wheels that
differ by a permutation lies in the subspace with an insulated
vertex. Thus the Alexander polynomial will not see the difference.
In fact the only thing that matters is the orientation of the wheel.
Therefore, we consider a clasper surgery on the unknot by a wheel
with $n$ legs attached by the simplest permutation. Pictured below
is the case of $n=8$.

\centerline{\includegraphics[width=2in]{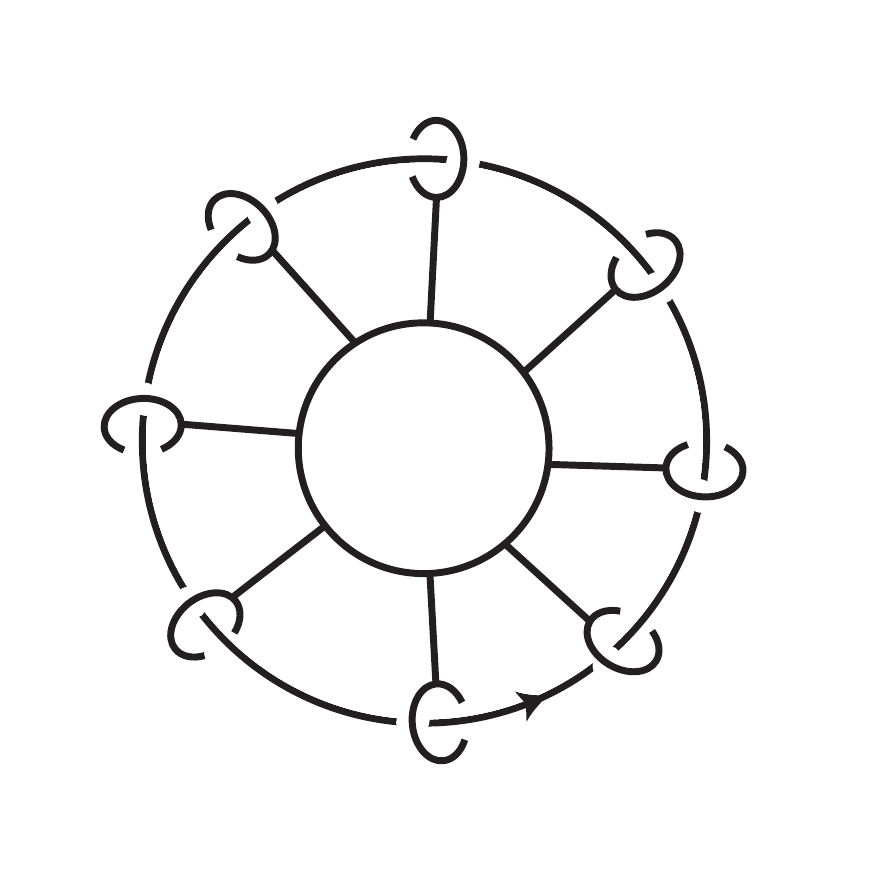}}

Let $\knot_n$ be the knot produced by this clasper surgery. The knot
$\knot_8$ is pictured in Figure~\ref{thisisaref}.
\begin{figure}
\begin{center}
\includegraphics[width=.4\linewidth]{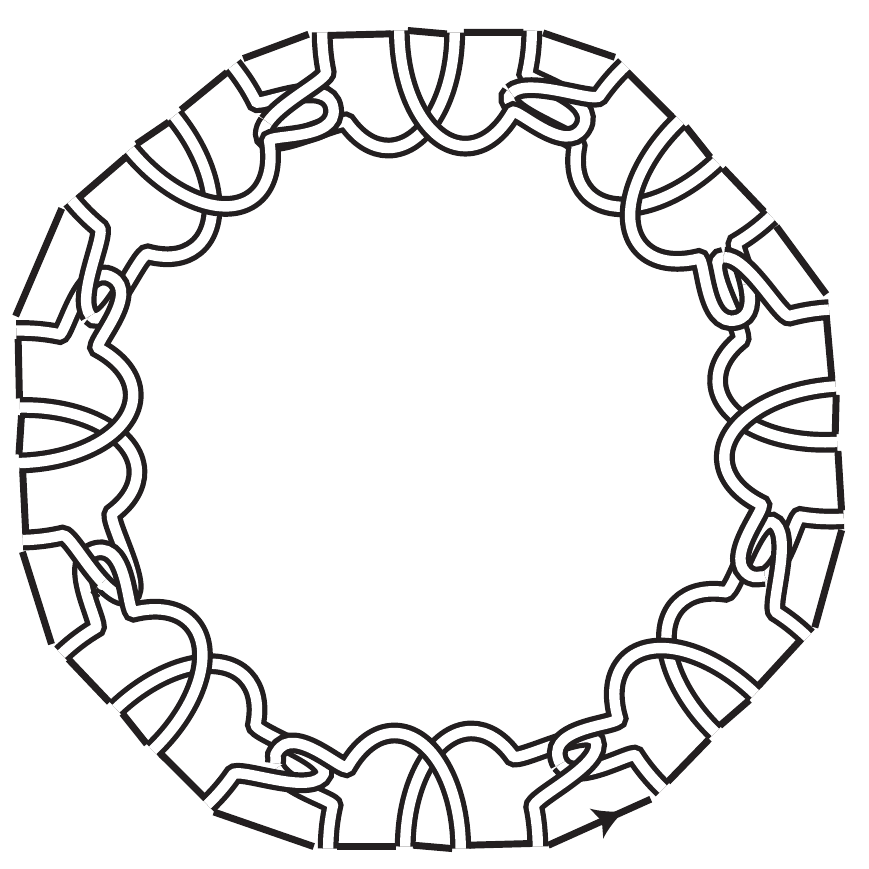}
\caption{The knot $k_8$. Note the evident Seifert
surface.}
\label{thisisaref}
\end{center}
\end{figure}

We prove that
\begin{itemize}
\item[1.] $pc_{2n}(\knot_{2n})=2$.
\item[2.] ${pc}_{2n}(\knot_{2n-1})\mod{2}$ is nonzero in $\mathbb Z/2$.
\end{itemize}
Because $\knot_{2n-1}$ and $\knot_{2n}$ both share finite-type invariants of order
$\leq(2n-2)$ with the unknot, it follows that $pc_{2n}=c_{2n}$
on these knots. The knot $\knot_n$ has an evident genus $n$ Seifert
surface with the following Seifert matrix:
$$\Theta=\left(\begin{array}{cc}0&M\\
M^T+I&0\end{array}\right),$$ where $M$ is the $n\times n$ matrix
$$M=\left(\begin{array}{ccccc}
0&1&0&\cdots &0\\
0&0&1&\cdots&0\\
\vdots&\ddots&\ddots&\ddots&\vdots\\
0&\cdots&\cdots&0&1\\
1&0&\cdots&0&0
                         \end{array}\right)$$
Now the Alexander polynomial can be calculated via the following
formula
$$A(\knot_n)=\det(t^{\frac{1}{2}}\Theta-t^{-\frac{1}{2}}\Theta^T)$$
which yields the result
$$A(\knot_n)=(1+(1-t)^n)(1+(1-t^{-1})^n).$$
See \cite{gl}.
The Conway polynomial $C(z)$ is obtained via the substitution
$z=t^{\frac{1}{2}}-t^{-\frac{1}{2}}$. Thus we have
$$A(\knot_n)=1+(t^{-n/2}+(-1)^nt^{n/2})z^n+(-1)^nz^{2n}.$$

Now, when $n=2k$, we get $1+(t^{-k}+t^k)z^{2k}+z^{4k}.$ Hence it
suffices to show that $t^{-k}+t^k$ is a function of $z$ with $2$ as
the leading term. First, this is function of $z$ because it is
symmetric with respect to $t$ and $t^{-1}$. To get the leading term,
plug in $z=0$ which is the same as $t=1$.

For $n=2k+1$, we get $1+(t^{-k-1/2}-t^{k+1/2})z^{2k+1}-z^{4k+2}$.
Thus it suffices to show that
$\phi(z)=t^{-k-1/2}-t^{k+1/2}=-(2k+1)z+\text{higher order terms}$.
Again, $t^{-k-1/2}-t^{k+1/2}$ is a function of $z$. But
\begin{align*}
t^{-k-1/2}-t^{k+1/2}&=(t^{-1/2}-t^{1/2})(t^{-k}+t^{-k+1}+\ldots+t^{k-1}+t^{k})\\
&=-z(2k+1+h.o.t.)
\end{align*}
This completes the proof.
\end{proof}

\begin{proof}[Proof of Proposition~\ref{nosystem}]
The fact that  $pc_{2n}$ is a Vassiliev invariant of order $2n$ has
already been observed. We now show that $pc_{2n}\mod 2$ is Vassiliev
of order $2n-1$.

To do this, it suffices to show that $pc_{2n}$ changes by an even
number after performing a degree $2n$ simple clasper surgery
(Goussarov \cite{gouss}, Habiro \cite{habiro}). The difference
between $pc_{2n}$ before and after the surgery is given by the value
of $pc_{2n}$ on the diagram in $\A^I_{2n}$ coming from the clasper.
Thus it suffices to show that the weight system for $pc_{2n}$ is
even on $\A^I_{2n}$. But this follows from Lemma~\ref{nosystem2}
since wheels generate $\A^I_{2n}$.

Lemma~\ref{nosystem2} also implies that the weight system for
$pc_{2n}\mod 2$ is given by $\ws$ and that $pc_{2n}$ induces the
weight system $2\hat{\ws}$.
\end{proof}

\subsection{Degree one}
Now we consider the question of which primitive Vassiliev invariants
are also $\DDelta_1$-invariants.

Let $\A^{2iv}$ be the submodule generated by diagrams with two
nonadjacent insulated vertices.
\begin{prop}
Let $f\colon \emph{\K}\to A$ be a Vassiliev invariant of order $n$
which is also a $\DDelta_1$-invariant. Then its weight system
$\hat{f}$, vanishes on $\A^{2iv}_n$.
\end{prop}
\begin{proof}
Proceed as in the proof of Proposition~\ref{iv}. If the two
insulated vertices are on different clasper components
$C_{i_o},C_{i_1}$, then the alternating sum is a sum of terms which
are themselves alternating sums of doing the $C_{i_o},C_{i_1}$
clasper surgeries. Since $C_{i_o}$, and $C_{i_1}$ are each sums of
$\DDelta$ moves, then $f$ must vanish on these alternating sums.

The argument from Proposition~\ref{iv} also works when the
two insulated vertices $v_0,v_1$  are on the same component,
provided that this pair of vertices is \emph{good}. Again
\emph{good} means that the complement of the two vertices, including
the Wilson loop, is connected.

We will show that every diagram with $2$ insulated vertices either
has a good pair, or is equivalent, modulo IHX, to one which has a
good pair.

The first step is that we can assume that our diagram has no
separating edges. In the argument from Proposition~\ref{iv} we
pushed the vertex into a trivalent graph $G$ which was a part of the
graph that has no separating edges, and is connected to the rest of
the diagram by a single separating edge. We argued $G$ had to have
at least one vertex, else the graph would have a loop edge, which is
zero modulo IHX. Indeed any symmetrical graph $G$ will have this
property, implying that $G$ has to have at least $6$ vertices. Hence
$G$ will have at least two vertices  $v_0$ and $v_1$ not joined by
an edge, which we can take to be our good pair.

So we may assume that the diagram minus the two vertices has exactly
two components. Let $H$ be the one not containing the Wilson loop.
Assume, without loss of generality, that the two vertices cut off an
$H$ of minimal size. Then every (trivalent) vertex of $H$ is
adjacent to both of $v_0$ or $v_1$. This leaves two possibilities:
$H$ is a ``Y" or $H$ is an ``H".

These two possibilities are pictured below:

\centerline {\includegraphics[width=.8\linewidth]{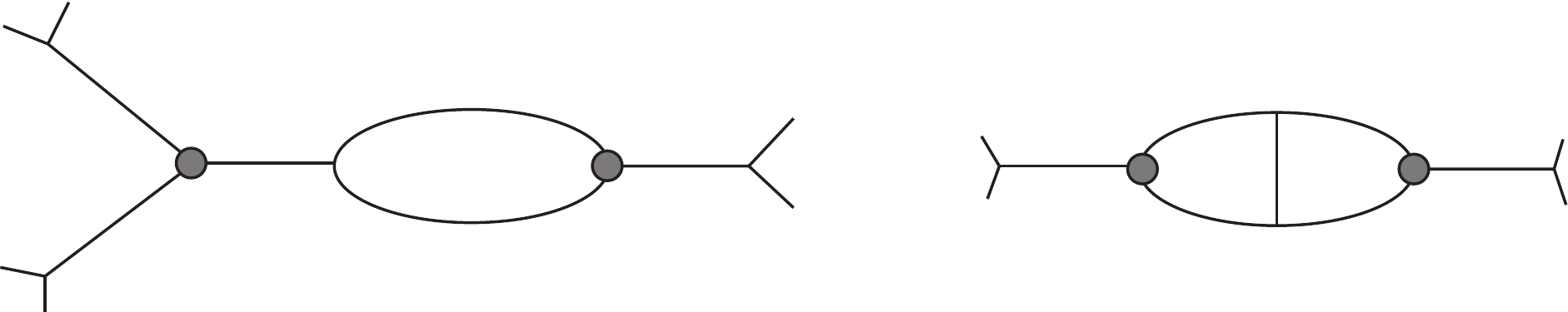}}

This first possibility means that we have a diagram $D$, with an
insulated vertex and a bubble inserted at an edge incident to that
vertex.The second possibility means that we have a diagram with an
edge not adjacent to the Wilson loop into which a theta graph has
been inserted.

Let us take care of the first possibility. It is known that
modulo IHX a bubble inserted in any edge is the same as inserting
it in any other edge (See for example \cite{ng}).
Hence insert a bubble on an edge incident to the Wilson loop. One
of the vertices of this bubble, $v^*$, is insulated. If this does
not form a good pair with any other insulated vertex, find another
insulated vertex which, together with $v^*$, separates off a minimal
graph $H^*$. If $H^*$ is an ``H", we are in the case of possibility
2. Otherwise we have a graph with two bubbles. Find an edge incident
to the Wilson loop and insert the bubbles on each of the two edges
joining it at a trivalent vertex. The two endpoints of the bubbles
near this trivalent vertex form a good pair, as pictured below:

\centerline{\includegraphics[width=.3\linewidth]{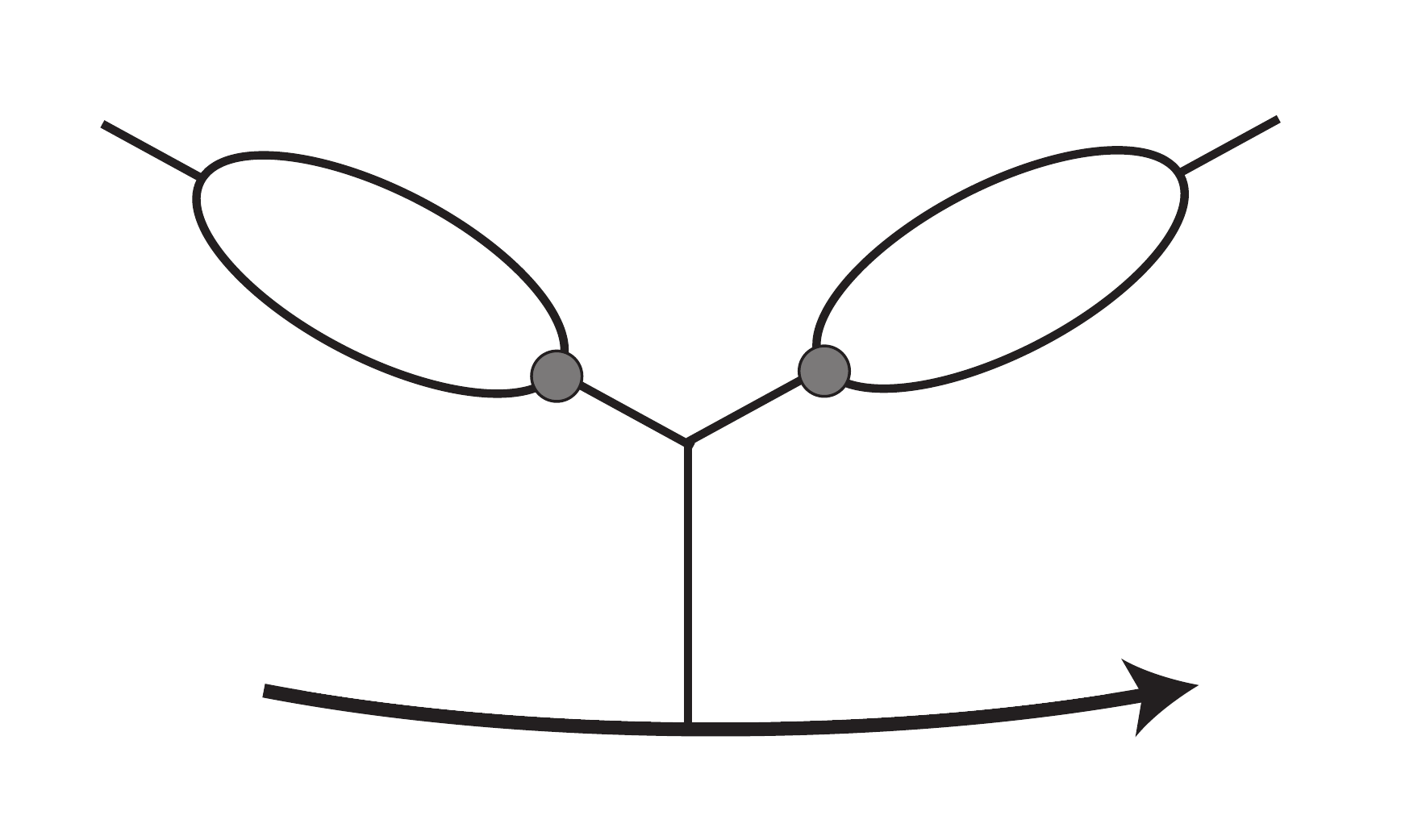}}

Now we consider the second possibility. The theta we inserted on an
edge is formed by first inserting a bubble in an edge and then
inserting a bubble in one of the two created vertices. By \cite{ng},
inserting a bubble at a vertex is independent of the
vertex.

Find a place where the diagram hits the knot. Let $v$ be the
adjacent trivalent vertex.  Insert a bubble at a vertex adjacent to
$v$ and insert a bubble in the opposite edge connecting into $v$.
The two vertices of the bubbles closest to $v$ are a good pair:

\centerline{\includegraphics[width=.3\linewidth]{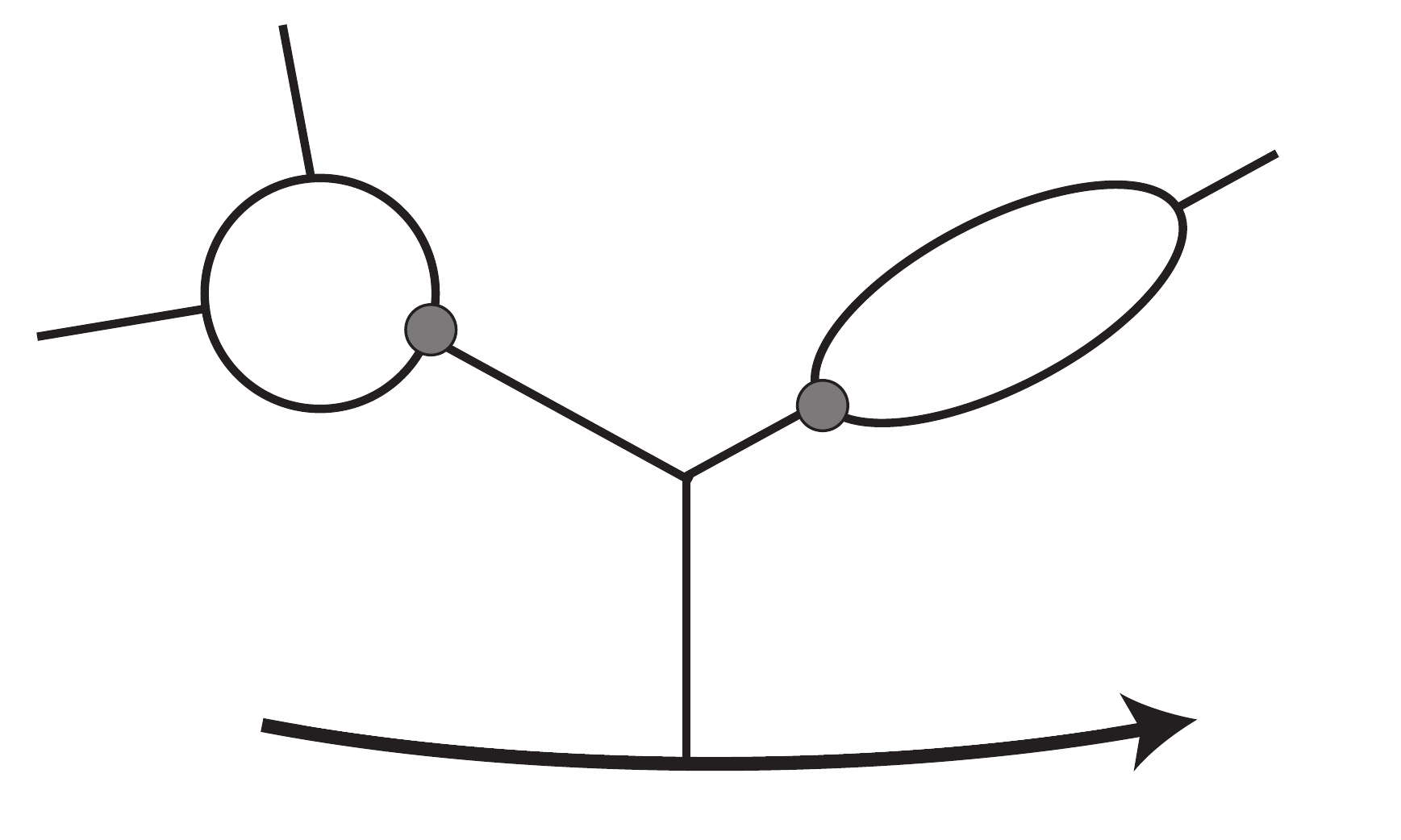}}

\end{proof}

We also have something of a converse.

\begin{prop}\label{conv}
The inclusion $\F^{\DDelta}_{k}\subset\F^{V}_{2k}$ holds.
Indeed,
any alternating sum of $k$ doubled-delta moves corresponds to a diagram in $\A_{2k}$
with $k$ insulated vertices, no two of which are joined by
an edge.
\end{prop}
\begin{proof}
A fact which we have been using repeatedly is that
any simple clasper surgery of degree $n$ will produce a diagram in $\A_n$ which is essentially gotten by interpreting the clasper as the internal graph and the knot as the Wilson loop.
(See \cite{ct2,habiro}.)
Using Habiro's zip construction, any capped clasper will produce a signed sum of diagrams where each leaf is turned into a univalent vertex which can attach to any of the arcs intersecting the leaf's cap. The sign is determined by the orientation of the cap and of the arc.

Hence each doubled-delta move produces a sum of ``Y"'s attached to the
Wilson loop in $8$ possible ways, by distributing the endpoints of
the $Y$ to two antiparallel arcs in the doubled-delta move. One can
visualize this as a $Y$ with two dashed edges, called \emph{tines}, emanating from each
univalent vertex and connecting to the appropriate places on the
Wilson loop.
The two tines have opposite sign, corresponding
to the antiparallelity. Each pair of tines is called a \emph{fork}. So every doubled-delta move produces a Y with a fork at each univalent vertex.

We simplify a diagram containing a Y with forks by sliding one tine of a fork
along the knot until it is adjacent to the other tine. The
resulting diagram is zero because it is a difference of two identical diagrams (with forks).
Now such a slide can be thought
of as a sum of diagrams, one for every edge attaching to the Wilson loop between
the fork's tines, where the fork is placed so as to straddle
each of these intervening edges. However, by the STU relation, such a straddling fork is the
same as attaching the corresponding $Y$ endpoint to the intervening
edge. Doing this for the other two forks in the $Y$ we get a sum of
diagrams where the endpoints of the $Y$ are connected to other
edges. So the Y's central vertex is insulated. Do this for all the
other $Y$'s sequentially. In the end we get $k$ insulated vertices.
Furthermore, no two will be joined by an edge.
\end{proof}

For the next argument, let $Z\colon\K\to \A\otimes\Q$ be the
Kontsevich integral, and let $\A^{\ell iv}$ be the submodule of $\A\otimes\Q$ generated
by diagrams with $\ell$ insulated vertices, no two of which share an
edge. (When $\ell=2$ this is a (minor) abuse of notation.)

\begin{prop}\label{whozwhatzit}
Let $C_1,\ldots, C_\ell$ be disjoint $ncY$-claspers on a knot $\knot$.
Then
$$Z[\knot;C_1,\ldots,C_\ell]\in\A^{\ell iv}.$$
\end{prop}
\begin{proof}
Let $T$ be a tangle comprised of three straight vertical strands, and let $C$ be a braid-like
Y-clasper whose leaves are meridians to the three strands.
Now let $T^d$ be the same tangle, but with each strand doubled into
a pair of antiparallel strands.
Assume (in the same spirit as the
Aarhus papers) that these antiparallel strands are infinitesimally
close.

{\bf Fact:} Let $\Xi=Z(T^d_C-T^d)$. Then $\Xi=d(Z(T_C-T))$, where $d$ is the
map that takes the endpoint of a diagram on $T$ and sums over
distributing it to the two new components. Since one of the strands
is oriented oppositely, it gets a minus sign.

The proof of this fact follows directly from Kontsevich's original
definition of the Kontsevich integral.

Note that the components of $\Xi$ each contain a trivalent
vertex because $C$ does not alter linking numbers.

Now $Z[\knot;C_1,\ldots,C_\ell]$ can be calculated as follows. Fix a height function and isotop the claspers $C_i$ so that they are surrounded by balls which look like the model tangle
$(T^d,C)$ above. Then the Kontsevich integral can be calculated locally and then the resulting diagrams glued together to give the Kontsevich integral of the entire knot.
In this way, we see that the Kontsevich integral will contain $\ell$ copies of $\Xi$, one for each of the balls surrounding each $C_i$.
Now each copy of $\Xi$ is a ``diagram
with forks" as in the proof of Proposition~\ref{conv}, where the
diagram contains at least one vertex.This gives rise to $k$
insulated vertices as in that proof.
\end{proof}

In order to identify $\A^I_n/\A^{2iv}_n$,
we first attack the question of what happens rationally. Let
$\B^I_n$ be the rational vector space spanned by connected,
vertex-oriented, unitrivalent graphs with $2n$ vertices, modulo the
antisymmetry and IHX relations. Let $\B^{2iv}_n$ be the subspace
spanned by graphs which have two vertices, neither of which is
joined by an edge to a univalent vertex, and which don't share an
edge with each other.

\begin{prop}
$\A^I_n/\A^{2iv}_n\otimes\Q\cong\B^I_n/\B^{2iv}_n$
\end{prop}
\begin{proof}
According to \cite{bnt}, there is an algebra isomorphism
$\chi\circ\partial_\Omega\colon \B\to \A\otimes\Q$, where $\B$ is
the space of not-necessarily-connected unitrivalent graphs, modulo
$IHX$ and $AS$. Let $\B^{2iv}$ be the subspace generated by graphs
which have two trivalent vertices not connected by an edge, neither
of which is connected by an edge to a univalent vertex.

We argue that
$\chi\circ\partial_\Omega(\B^{2iv})=\A^{2iv}\otimes\Q$. Clearly
$\chi\circ\partial_\Omega(\B^{2iv})\subset\A^{2iv}\otimes\Q,$ by
the definitions of $\chi$ and $\partial_\Omega$. For the reverse
direction, note that $\chi^{-1}=\sigma$ (See \cite{bn}), is defined
by an iterative process in which diagrams keep getting more things
tacked onto them. In particular, insulated vertices will stay
insulated. Thus $\sigma(\A^{2iv}\otimes\Q)\subset\B^{2iv}$. Now we
wish to show that $\partial_\Omega^{-1}(\B^{2iv})\subset \B^{2iv}$.
So let $A$ be such that $\partial_\Omega(A)\in\B^{2iv}$. Suppose,
toward a contradiction, that $A\not\in\B^{2iv}$.
Consider the smallest degree where $A$ has terms not in $\B^{2iv}$.
 Let $s$ be the sum
of terms in $A$ which are of this degree. Then $\partial_\Omega(s)=s+h.o.t.$ Hence
the terms in
$s$ which are not in $\B^{2iv}$ cannot
cancel with anything else, which is the desired contradiction.

Now $\chi\circ\partial_\Omega$ will send disconnected graphs to
separated graphs in $\A\otimes\Q$, thus inducing an isomorphism
$\B^I_n/\B^{2iv}_n\to\A^I_n/\A^{2iv}_n\otimes\Q$.
\end{proof}

\begin{prop}
There is an isomorphism $\B^I_n/\B^{2iv}_n\cong\Q$.
\end{prop}
\begin{proof}
Elements of $\B^I_n$ consist of a trivalent core graph with
edges that have a univalent vertex
attached to some of the edges. These edges with univalent vertices are called \emph{hairs}.

If the core graph is not a circle, a
theta, or two loops joined by an edge (barbell graph), then we have an element of
$\B^{2iv}$. However, the IHX relation can be used to show that
a barbell graph with hairs is always zero. (If there is no hair on the non-loop edge, apply an IHX relation to that edge. If there is at least one hair on the non-loop edge, apply an
IHX relation to distribute a hair close to a loop to the loop in two canceling ways.)

Thus $\B^I_n/\B^{2iv}_n$ is
generated by a wheel together with graphs that have theta as a core
graph. When $n$ is even, this latter type of graph has an
orientation reversing symmetry, so we are left with just the wheel.
So we get $\Q$ in that case. On the other hand, if $n$ is odd, the
wheel has an orientation reversing symmetry, and disappears. Thus we
need to analyze the space of theta graphs with hairs attached to the
edges, modulo $\B^{2iv}_n$.
These can be characterized as unordered
triplets $(a,b,c)$, where $a,b,c$ represent the number of hairs
attached to each edge. Hence $a+b+c=n-1$.
There are two types of IHX relation. The first kind can occur when $a,b$ or $c$ is zero, and alters the structure of the
underlying core graph. However this type of IHX relation merely asserts that a certain barbell graph is zero; the two terms with a theta as core graph cancel.

The second kind of IHX relation turns
into the following relation
$(a+1,b,c)+(a,b+1,c)+(a,b,c+1)=0$
for any $a+b+c=n-2$. Also $(a,b,c)\in\B^{2iv}_n$ if
all of $a,b,c$ are nonzero. Thus the only relations that don't lie
entirely in $\B^{2iv}_n$ are for triples $(a,b,c)$,
where $a=0$. If $b=0$, as well, the ensuing relation is
just $(0,0,n-1)+2(0,1,n-2)=0.$ Otherwise, modulo $\B^{2iv}$, the
relation is just $(0,a,b)=(0,a-1,b+1)$, where $b>a\geq 2$. Thus we
see that
 $\B^I_n/\B^{2iv}_n\cong \Q$, generated by $(0,1,n-2)$.
\end{proof}

This allows us to do a partial calculation over $\Z$.

\begin{prop}
We have
\begin{align*}
\Z\oplus\Z_2&\twoheadrightarrow\A^I_{2n}/\A^{2iv}_{2n}\\
\Z&\cong\A^I_{2n-1}/\A^{2iv}_{2n-1}
\end{align*}
\end{prop}
\begin{proof}
As always, we use that $\A^I_n$ is generated by wheels attached to
the Wilson loop by some permutation. Let $Y\in\A^I_n$ be a theta
graph with a single hair on one edge, no hairs on the second edge
and $n-2$ hairs on the third edge, attached to the Wilson loop by some
permutation.
 Modulo $\A^{2iv}_n$, the order in
which the legs is attached to the Wilson loop is irrelevant. So,
abusing notation, let these all be called $Y$. Then any two ways of
attaching a wheel to the Wilson loop will differ, modulo
$\A^{2iv}_n$ by multiples of $Y$. Let $X$ be the wheel attached by
the identity permutation. Then we have just shown that $X$ and $Y$
generate.

Now in the case when $n$ is even, $Y$ has an odd number of trivalent
internal vertices. Hence, flipping $Y$ over (and reordering the legs
accordingly) implies that $Y=-Y$, or that $2Y=0$. On the other hand,
we know $X$ has infinite order, since it corresponds to $pc_{n}/2$.
This establishes the first statement above.

Suppose, then, that $n$ is odd. Flip $X$ over, and rearrange the
legs to be the identity permutation again. Counting the number of
$Y$'s we get as error terms, we get the equation $X=-X+(2n-3)Y$,
implying $2X=(2n-3)Y$. Thus $\A^I_{2n-1}/\A^{2iv}_{2n-1}$ is a
quotient of $\Z\oplus\Z/\langle2,2n-3\rangle\cong \Z$. On the other
hand, we know that tensoring with $\Q$ gives us $\Q$, so that in
fact no further relations are possible.
 \end{proof}

This proves the first part of Theorem \ref{dd1} of the introduction.

\begin{cor}
There exist Vassiliev invariants of order $2n-1$
$v_{2n-1}\colon\emph{\K}\to\Q$ which are $\DDelta_1$-invariants.
\end{cor}
\begin{proof}
Let $W_{2n-1}\colon\A^I_{2n-1}\to \Z$ be a weight system that
vanishes on $\A^{2iv}_{2n-1}$ which detects
$\A^{I}_{2n-1}/\A^{2iv}_{2n-1}\cong \Z$.  Let
$v_{2n-1}=W_{2n-1}\circ Z_{2n-1}$. Clearly $v_{2n-1}$ is a Vassiliev invariant
of order $2n-1$. Let $C_1,C_2$ be two $\DDelta$-moves on a knot
$\knot$. Then $ Z_{2n-1}[\knot;C_1,C_2]\in \A^{2iv}_{2n-1}$ by
Proposition~\ref{whozwhatzit}.
\end{proof}

\begin{prop}
There is an isomorphism $\A_4^I/\A^{2iv}\cong \Z\oplus \Z/2$.
\end{prop}

\begin{proof}
One uses the generators $X$ and $Y$, where $X$ is a wheel, say,
attached by the identity permutation, and $Y$ is a tree of degree
$4$ with no crossings. By the usual sliding argument, any diagram
can be realized, modulo separated diagrams, as a sum of diagrams
with at least one loop. Moreover, this will preserve any insulated
vertices, and possibly create new ones. So we list all nonzero
graphs that have a core equal to a trivalent graph, $G$. The
possible $G$'s are:
\begin{center}
\includegraphics[width=.5\linewidth]{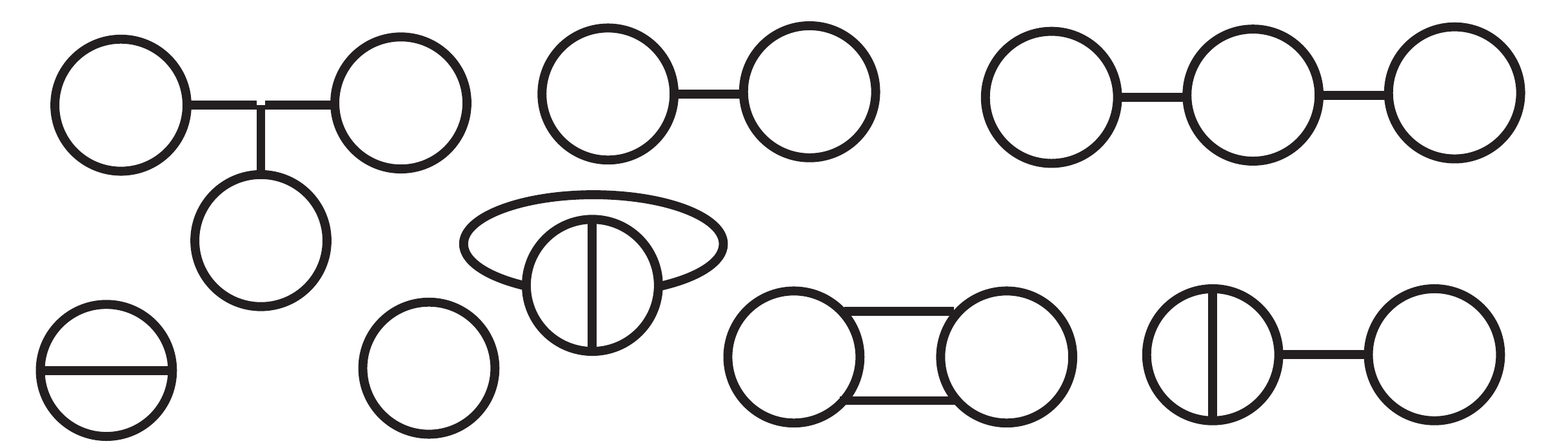}
\end{center}
Add forests to these cores in all ways, and enumerate those that
have two insulated vertices. All of them are an even multiple of
$Y$, and indeed the graph which is a tetrahedron with one hair on
two adjacent edges, is equal to $\pm 2Y$.
\end{proof}

\begin{cor}\label{Sequivintro2}
Any integer-valued finite-type invariant of degree $4$ is, modulo
two, a doubled-delta invariant of degree $1$.
\end{cor}
\begin{proof}
According to Proposition~\ref{conv}, it suffices to notice that any
degree $4$ finite-type invariant vanishes on $\A^{2iv}$ modulo $2$.
\end{proof}

This establishes the second part of Theorem \ref{dd1} of the introduction.

   %%%%%%%%%%%%%%%%%%%%%%%%%%%%%%%%%%%%%%%%%%%%%%%%%%%%%%%%%%
   % ¡mambo!                                                %
   %                                                        %
   %      %%%%%%%%%%%%%%    %%%%%    %%%%%       %%%%%      %
   %       %%%        %      %%%      %%%         %%%       %
   %       %%%               %%%      %%%%        %%%       %
   %       %%%               %%%      %%% %       %%%       %
   %       %%%               %%%      %%%  %      %%%       %
   %       %%%               %%%      %%%   %     %%%       %
   %       %%%    %          %%%      %%%    %    %%%       %
   %       %%%%%%%%          %%%      %%%     %   %%%       %
   %       %%%    %          %%%      %%%      %  %%%       %
   %       %%%               %%%      %%%       % %%%       %
   %       %%%               %%%      %%%        %%%%       %
   %       %%%               %%%      %%%         %%%       %
   %      %%%%%             %%%%%    %%%%%       %%%%%      %
   %                                                        %
   %%%%%%%%%%%%%%%%%%%%%%%%%%%%%%%%%%%%%%%%%%%%%%%%%%%%%%%%%%

\end{document}